%% file: manuscript.tex
\documentclass[onecolumn,
]{mpi2015-cscpreprint}

\usepackage{amsthm,amsmath,amssymb,amsfonts}
\usepackage{mathtools}

\usepackage{xcolor}
\usepackage{graphics}
\usepackage{graphicx}
\usepackage{hyperref}
\usepackage{float}
\usepackage{placeins}
\usepackage{booktabs, multirow}
\usepackage{url}
\usepackage{verbatim,upquote}
\usepackage[section]{algorithm}

\usepackage{multirow}
\usepackage{threeparttable}

\usepackage{subcaption}
\numberwithin{figure}{section}

\usepackage{listings}
\usepackage{upquote}
\usepackage{tabularx}
\usepackage{pgfplots}
\pgfplotsset{compat=1.18}
\usepgfplotslibrary{groupplots}

\newenvironment{codefragment}[1][htb]{
	\begin{algorithm2e}[#1]}{\end{algorithm2e}}

\usepackage[algo2e,ruled,linesnumbered,noline,algosection]{algorithm2e}
\DontPrintSemicolon
\makeatletter
\let\oldnl\nl \newcommand{\nonl}{\renewcommand{\nl}{\let\nl\oldnl}} \renewcommand{\SetKwProg}[4]{\algocf@newcmdside@koif{#1}{\KwSty{#2}\ifArgumentEmpty{#2}\relax{\ 
		}\ProgSty{##2}\KwSty{#3}\ifArgumentEmpty{##1}\relax{ 
			##1}\\\SetAlgoVlined##3\SetAlgoNoLine{#4}{##4}}\algocf@newcommand{l#1}{\@ifstar{\csname algocf@l#1star\endcsname}{\csname 
			algocf@l#1\endcsname}}\algocf@newcmdside{algocf@l#1}{3}{\KwSty{#2} 
		\ProgSty{##2}\KwSty{#3}\algocf@bgroupcode\ 
		##3\algocf@egroupcode\@endalgocfline\ifArgumentEmpty{##1}\relax\ 
		{##1}\strut\par}\algocf@newcmdside{algocf@l#1star}{3}{\KwSty{#2} 
		\ProgSty{##2}\KwSty{#3}\algocf@bgroupcode\ ##3\algocf@egroupcode}}
\newcommand{\algorithmfootnote}[2][\footnotesize]{\let\old@algocf@finish\@algocf@finish \def\@algocf@finish{\old@algocf@finish \leavevmode\rlap{\begin{minipage}{\linewidth}
				#1#2
		\end{minipage}}}}
\SetKwInOut{Parameter}{Parameter}

\newcommand{\algorithmstyle}[1]{\renewcommand{\algocf@style}{#1}}
\makeatother
\SetKwProg{Fn}{function}{}{}
\SetProgSty{textnormal}
\SetKwProg{subFn}{subfunction}{}{}

\SetKwComment{Comment}{$\triangleright$\ }{}
\SetCommentSty{itshape}
\newcommand{\myhrulefill}{\leavevmode\leaders\hrule height 0.7ex depth 
	\dimexpr0.4pt-0.7ex\hfill\kern0pt}
\newcommand{\lineseparator}{\nonl\myhrulefill\\}

\graphicspath{{figs/}}

\newtheorem{theorem}{Theorem}
\newtheorem*{theorem*}{Theorem}

\numberwithin{equation}{section}
\numberwithin{table}{section}
\numberwithin{figure}{section}
\numberwithin{theorem}{section}

\definecolor{hotpink}{rgb}{0.9,0,0.5}
\hypersetup{colorlinks,urlcolor=blue,citecolor=hotpink,linkcolor=blue}

\def\eu{\ensuremath{\mathrm{e}}}

\def\>{\mskip\medmuskip}

\DeclarePairedDelimiter\norm{\lVert}{\rVert}
\newcommand{\normi}[1]{\ensuremath{\norm{#1}_1}}
\newcommand{\normF}[1]{\ensuremath{\norm{#1}_F}}

\def\normt#1{\|#1\|_{2}}

\newcommand{\mPi}{\ensuremath{\mathit{\Pi}}}

\newcommand{\R}{\ensuremath{\mathbb{R}}}

\newcommand{\nbyn}{\ensuremath{n\times n}}

\newcommand{\ldlt}{\ensuremath{LDL^T}}
\renewcommand{\L}{\ensuremath{\mathcal{L}}}

\DeclareMathOperator{\diag}{diag}
\DeclareMathOperator{\cond}{cond}

\DeclareMathOperator{\rank}{rank}

\DeclareMathOperator{\sign}{sign}
\DeclareMathOperator{\size}{size}
\DeclareMathOperator{\blkdiag}{blkdiag}

\DeclareMathOperator{\vvec}{vec}

\DeclareMathOperator{\fl}{\operatorname{f\kern.2ptl}}

\newcommand{\us}{u_s}

\newcommand{\wh}{\widehat}
\newcommand{\wt}{\widetilde}
\newcommand{\gammat}{\ensuremath{\tilde{\gamma}}}

\mathcode`@="8000 {\catcode`\@=\active\gdef@{\mkern1mu}}

\newcounter{mylineno}
\makeatletter
\let\oldtabcr\@tabcr

\def\mynewline{\refstepcounter{mylineno}\llap{\footnotesize\arabic{mylineno}\hspace{5pt}}}

\gdef\@tabcr{\@stopline \@ifstar{\penalty \@M \@xtabcr}\@xtabcr\mynewline}

\begin{document}

\title{Mixed-precision iterative refinement for low-rank Lyapunov equations\thanks{
		Version of October 30, 2026. Open Access funding provided by the Max Planck Society within the Projekt DEAL agreement.
		}
}
 
\author[$\dagger$]{Peter Benner}
\affil[$\dagger$]{Max Planck Institute for Dynamics of Complex Technical Systems, Magdeburg, Germany\authorcr
	\email{benner@mpi-magdeburg.mpg.de}, \orcid{0000-0003-3362-4103}}

\author[$\dagger$,$\ast$]{Xiaobo Liu}
\affil[$\dagger$]{Max Planck Institute for Dynamics of Complex Technical Systems, Magdeburg, Germany\authorcr
	\email{xliu@mpi-magdeburg.mpg.de}, \orcid{0000-0001-8470-8388}}

\affil[$\ast$]{Corresponding author.}

\shortauthor{P. Benner and X. Liu}
  

  
\abstract{
We develop a mixed-precision iterative refinement framework for solving low-rank Lyapunov matrix equations $AX + XA^T + W =0$, where $W=LL^T$ or $W=LSL^T$. Via rounding error analysis of the algorithms we derive sufficient conditions for the attainable normwise residuals in different precision settings and show how the algorithmic parameters should be chosen. These conditions are independent of the choice of inner solver, provided that the prescribed residual accuracy is attained in the inner solves. 
Using the sign-function Newton iteration as the solver, we demonstrate that reduced precisions, such as half precision with unit roundoff $u_s$, can be used for Lyapunov equations with condition number of order $u_s^{-1}$ without compromising the attainable solution quality. 
The numerical experiments validate our error analysis, and we establish cost models to quantify the expected speedup for different low-to-high precision timing ratios. 
This work provides an algorithmic framework towards exploiting native low-precision hardware to accelerate low-rank Lyapunov solvers without sacrificing accuracy.
}

\novelty{
We develop a mixed-precision IR framework for the factored solution of low-rank Lyapunov equations, in the formulation of either Cholesky-type or $\ldlt$-type.
We provide new rounding error analysis, which indicates how to set the precisions and choose the algorithmic parameters within the IR framework. 
The experiments demonstrate the potential of exploiting reduced precisions, such as the half precision, to accelerate the solution of low-rank Lyapunov equations.
}

\maketitle

\section{Introduction}\label{sect.intro}
This paper studies the mixed-precision iterative refinement (IR) algorithm for low-rank continuous-time Lyapunov equations, which have the form
\begin{equation}\label{eq:lyap}
	AX + XA^T + W =0, \quad A, W\in\R^{\nbyn}.
\end{equation}
Such equations arise in the computation of controllability and observability Gramians of linear time-invariant dynamical systems~\cite{anto05,zdg96}. These Gramians play key role in system balancing~\cite{lhpw87} and balanced truncation model order reduction~\cite{moor81,pesi82}, including methods for large-scale
data-sparse systems~\cite{babe08}.
Here $A$ is Hurwitz (asymptotically stable), and $W$ is symmetric positive semidefinite with rank $m\ll n$, admitting a Cholesky- or $LDL^T$-type low-rank factorization 
\begin{equation}\label{eq:lyap-W-fac}
	W=LL^T \quad \mbox{or} \quad W=LSL^T, \quad
	L\in\R^{n\times m},
\end{equation}
where $S$ is symmetric positive semidefinite. 
The stability of the coefficient matrix $A$ implies a unique symmetric positive semidefinite solution, given by the integral representation~\cite{hamm82}
\begin{equation}\label{eq:lyap-X-intgral}
 X=\int_0^\infty e^{tA}W e^{tA^T}\, \mathrm{d}t,
\end{equation}
which shows that $X$ need not have low rank, despite the low rank of $W$. 
In particular, if $(A,LS^{1/2})$ is controllable, then $X$ is positive definite and has full rank~\cite[Prop.~3.1]{asz02}.
The solution $X$ can have low rank when the range of $W$ is contained in an $A$-invariant subspace of much smaller dimension, say, $r\ll n$.
In this case, every integrand in~\eqref{eq:lyap-X-intgral} then has its range contained in that subspace, and hence $\rank(X)\le r$. 

More generally, even a full-rank solution $X$ can have low \emph{numerical rank}, namely, it can be accurately approximated by a low-rank matrix. 
This can happen when the eigenvalues of the solution decay sufficiently rapidly, and  such decay can be quantified, for example, when $A$ is symmetric negative definite~\cite{penz00b}.
For such $A$, if its spectrum is contained in $[-\beta,-\alpha]$, where $0<\alpha<\beta$, and $0<\rank(W)\le m$, then one can bound the error in the best approximation $X_{mq}$ of rank at most $mq$~\cite[Thm.~2.1, Cor.~3.2]{beto17} by
\begin{equation*}
 \frac{\normt{X-X_{mq}}}{\normt{X}}
 \le 4\exp\left(-\frac{\pi^2 q}{\log(4\beta/\alpha)}\right),
 \qquad mq<n,\quad q\in\mathbb{N},\quad q\ge1.
\end{equation*}
This gives an exponentially decaying upper bound on the singular values of $X$ for fixed $m$ and $\beta/\alpha$.
For general nonnormal $A$, the decay of the singular values of $X$ cannot be predicted from the eigenvalues of $A$ alone, but it may also depend on the conditioning of its eigenvectors~\cite{asz02} and on its numerical range or pseudospectra~\cite{bes15}.

When $X$ can be well-approximated by a low-rank matrix, a fundamental and well-established strategy in the literature is to seek an approximate solution directly in factored form, 
\begin{equation}\label{eq:lyap-X-fac}
	X\approx ZZ^T \quad\mbox{or}\quad X\approx ZYZ^T,
	\qquad Z\in\R^{n\times r},\quad Y\in\R^{r\times r},\quad r\ll n.
\end{equation}
This idea of seeking factored solutions can be traced back to Hammarling~\cite{hamm82}, who exploited it for solving stable and non-negative definite Lyapunov equations.
Complex dynamical systems can lead to equation~\eqref{eq:lyap} with dimension of order $10^5$ or larger~\cite{bekq11}.
For these large-scale Lyapunov equations, iterative solvers that compute and update the low-rank factors in~\eqref{eq:lyap-X-fac} are typically preferred, since factorization-based methods become unduly expensive in terms of both computational overhead and memory storage.
Such numerical schemes include the method based on rational iterations for the matrix sign function~\cite{bequ99}, the cyclic low-rank Smith method~\cite{penz00a}, the low-rank alternating direction implicit (LR-ADI) method~\cite{liwh04,blp08}, and projection-type methods based on Krylov subspaces~\cite{jaka94},~\cite{jbri06},~\cite{simo07}; see~\cite{besa13,simo16} for surveys.

Modern hardware increasingly supports native precisions lower than the traditional IEEE binary64 (fp64) and binary32 (fp32) formats~\cite{ieee19}, and this has fostered the development of mixed-precision algorithms.
Utilizing the low precisions appropriately within numerical algorithms can accelerate computation, reduce data storage and communication, and improve energy efficiency on computational units, without sacrificing their accuracy and stability. 
We refer the reader to~\cite{hima22} for a survey on recent developments of mixed-precision algorithms in numerical linear algebra.

With conventional floating-point units, half precisions, including the IEEE binary16 format (fp16) and the bfloat16 format (bf16) by Google Brain,\footnote{https://research.google} can offer approximately $2\times$ or $4\times$ the arithmetic throughput of fp32 or fp64, respectively.
The advent of tensor-core accelerators has however pushed the theoretical peak throughput of fp16 matrix operations on the V100 to $8\times$ or $16\times$ that of fp32 or fp64, respectively~\cite{htd19}, and practical performance evaluation on this GPU has revealed that tensor cores can boost the \texttt{GEMM} (general matrix--matrix multiply) performance over fp32 CUDA-core implementations by up to $6\times$ when multiplying large matrices and $12\times$ when
multiplying small-size matrices in parallel~\cite{mclp18}. 
In particular, using fp16 tensor cores within a GMRES-based
mixed-precision IR scheme can provide up to $4\times$ speedup
over the fp64 LU-based \texttt{dgesv} solver on a V100
GPU~\cite{htdh18}, and an asymptotic $3\times$ speedup
over a full-fp64 LU-based direct solver on an A100 GPU~\cite{abcg24}, while delivering fp64 accuracy.
On more recent hardware, a single NVIDIA B200 GPU has a peak throughput of 2.25 PFLOP/s for dense fp16/bf16 Tensor Core matrix operations, compared with
75 TFLOP/s for fp32 and 37 TFLOP/s for fp64, corresponding to peak throughput ratios of about $30\times$ and $61\times$, respectively~\cite[p.~8]{nvid25}.
These hardware peak rates illustrate the huge potential of accelerating scientific computations using reduced precision.
Table~\ref{table:fp-precs} reports the key parameters of the four floating-point arithmetics considered in this work.

\begin{table}[t]
	\caption{Parameters of bf16, fp16, fp32, and fp64 arithmetic:
		significand precision $t$ (including the implicit bit),
		exponent-field width $e$, unit roundoff $u=2^{-t}$,
		smallest positive normal number $x_{\min}=2^{e_{\min}}$,
		and largest finite number
		$x_{\max}=(2-2^{1-t})2^{e_{\max}}$,
		where $e_{\min}=2-2^{e-1}$ and
		$e_{\max}=2^{e-1}-1$.
		Numerical values are given to three significant figures.}
	\label{table:fp-precs}
	\centering
	\setlength{\tabcolsep}{8pt}
	\begin{tabular}{l r r l l l}
		\toprule
		\qquad Format & $t$ & $e$ & \qquad$u$ & \quad$x_{\min}$ & \quad$x_{\max}$ \\
		\midrule
		bfloat16\phantom{i} (bf16) & 8 & 8
		& $3.91\times 10^{-3}$
		& $1.18\times 10^{-38}$
		& $3.39\times 10^{38}$ \\
		binary16 (fp16) & 11 & 5
		& $4.88\times 10^{-4}$
		& $6.10\times 10^{-5}$
		& $6.55\times 10^{4}$ \\
		binary32 (fp32) & 24 & 8
		& $5.96\times 10^{-8}$
		& $1.18\times 10^{-38}$
		& $3.40\times 10^{38}$ \\
		binary64 (fp64) & 53 & 11
		& $1.11\times 10^{-16}$
		& $2.23\times 10^{-308}$
		& $1.80\times 10^{308}$ \\
		\bottomrule
	\end{tabular}
\end{table}

Despite the wide use of mixed precision in the numerical linear algebra community, its potential has remained largely unexploited in approximating the solution of matrix equations.
Benner et~al.~\cite{bekq11} developed an algorithm for computing factored solution of the low-rank Lyapunov equation~\eqref{eq:lyap}, where the idea of IR in fp32 and fp64 was exploited, albeit not fully spelled out; 
their focus was more on the implementation with hybrid CPU--GPU platforms rather than algorithmic development.
More recently, based on a fixed-precision IR scheme for the  perturbed quasi-triangular Sylvester equation, a mixed-precision Schur-based method was devised for computing the full solution of the Sylvester matrix equation~\cite{dfhl25}.

In this paper, we develop a mixed-precision IR framework for solving low-rank Lyapunov matrix equations.
We provide a rounding error analysis of the algorithms to guide the choice of the precisions and algorithmic parameters. We examine the IR framework by using the sign function Newton iteration as the solver. 
We begin with the Cholesky-type and $\ldlt$-type mixed-precision IR frameworks in Section~\ref{sect.mp-ir}, followed by a rounding error analysis for the $\ldlt$-type formulation and its specialization to the Cholesky-type scheme.
In Section~\ref{sect.sign-func-Newton-iter} we discuss the use of the sign function Newton iteration within the IR framework.
Numerical experiments are presented in Section~\ref{sect.experiments} to verify the analysis and assess the quality of the computed solutions, together with estimates of the potential speedup of the mixed-precision IR framework based on floating-point operation (flop) counts. Conclusions are drawn in Section~\ref{sect.conclusions}.

Throughout the paper, we use MATLAB-style colon notation to denote index ranges; for example, $A(:, j)$ selects the entire $j$-th column, and $A(i, :)$ selects the entire $i$-th row.
The spectral radius of a square matrix $A$ is denoted by 
\[
\rho(A):=\max\bigl\{|\lambda| : \lambda\ \mbox{is an eigenvalue of}\ A\bigr\}.
\]
The $\diag(\cdot)$ operator creates a diagonal matrix from its input scalar elements, while $\blkdiag(\cdot)$ returns a block-diagonal matrix from its input matrices; and $\norm{\cdot}$ denotes any consistent operator norm.

\newcommand{\muk}[1][k]{\ensuremath{\mu_{#1}}}
\newcommand{\mukm}[1][k-1]{\ensuremath{\mu_{#1}}}
\newcommand{\chol}{\texttt{chol}} 

\section{Mixed-precision IR framework}\label{sect.mp-ir}
The low-rank Lyapunov equation~\eqref{eq:lyap} can be recast as the $n^2\times n^2$ Kronecker linear system
\begin{equation}\label{eq:lyap-ls}
 M\vvec(X)=w,\qquad M:=I_n\otimes A+A\otimes I_n,\quad w:=\vvec(-W),
\end{equation}
where $I_n$ denotes the identity matrix of order $n$, and $\vvec$ stacks the columns of an $m\times n$ matrix into a vector of length $mn$.
In theory, one can apply any linear system solver to the equivalent system~\eqref{eq:lyap-ls} for the solution of the Lyapunov equation~\eqref{eq:lyap}. This approach should nonetheless be avoided in practice, not only due to its prohibitively expensive storage requirements, but also because it is unclear how the low-rank structure can be exploited.
We therefore refine the solution factors directly, starting with the existing Cholesky-type IR scheme and then extending it to an $\ldlt$-type formulation. We carry out the rounding error analysis for the latter and conclude by specializing the analysis to the Cholesky-type scheme.
Throughout the rest of the paper,
\[
 \L\colon\R^{n\times n}\mapsto \R^{n\times n},\qquad \L(X)=AX+XA^T,
 \qquad \kappa_F(\L)=\normF{\L}\normF{\L^{-1}},
\]
where $\normF{\L}$ denotes the operator norm induced by the
Frobenius norm on $\R^{n\times n}$.

\subsection{Existing Cholesky-type IR}\label{sect:chol-ir}
The authors of~\cite{bekq11} design an IR scheme for the factored solution of the low-rank Lyapunov equation~\eqref{eq:lyap}, where the solver step is carried out in fp32 and the other steps are performed in the usual fp64 environment.
The major difference between this IR scheme and that for a linear system lies in the residual computation and solution update steps (see Lines~\ref{alg.line.ir3-lyap.resfac.chol} and~\ref{alg.line.ir3-lyap.solupt.chol} of Algorithm~\ref{alg.ir3-lyap-chol}): these steps involve more complex computational kernels, such as QR factorizations and spectral decompositions.

The residual of an approximated Cholesky-type factor $Z_i\in\R^{n\times c_i}$ of the solution to~\eqref{eq:lyap}, with $W=LL^T$, has the form
\begin{equation}\label{eq:lyap-res-chol}
 \mathcal R(Z_i):=AZ_iZ_i^T+Z_iZ_i^TA^T+LL^T
                =\L(Z_iZ_i^T)+LL^T.
\end{equation}
This residual can be indefinite due to the inexactness of the solver as well as the rounding errors. It can be rewritten as the product
\[
 \mathcal R(Z_i)=F_iP_iF_i^T,\qquad
 F_i=\begin{bmatrix}Z_i&AZ_i&L\end{bmatrix},\qquad
 P_i=\begin{bmatrix}0&I_{c_i}&0\\I_{c_i}&0&0\\0&0&I_m\end{bmatrix}.
\]
Thus the residual factorization can be performed via a thin QR factorization $F_i=U_iT_i$, followed by a spectral decomposition of the small kernel matrix $H_i=T_iP_iT_i^T$, without explicitly forming an $n\times n$ matrix.
The Cholesky-type IR scheme of~\cite{bekq11} splits this residual into positive semidefinite (PSD) and negative semidefinite (NSD) parts,
\begin{equation}\label{eq:res-fac}
 \mathcal R(Z_i)=L_i^+(L_i^+)^T-L_i^-(L_i^-)^T
                =:\mathcal R^+(Z_i)-\mathcal R^-(Z_i),
\end{equation}
and solves two PSD correction equations,
\begin{subequations}\label{eq:lyap-cor-chol}
\begin{align}
 AX_i^++X_i^+A^T+L_i^+(L_i^+)^T&=0,
 &X_i^+&=Z_i^+(Z_i^+)^T,\label{eq:lyap-cor-chol-plus}\\
 AX_i^-+X_i^-A^T+L_i^-(L_i^-)^T&=0,
 &X_i^-&=Z_i^-(Z_i^-)^T.\label{eq:lyap-cor-chol-minus}
\end{align}
\end{subequations}
Here $Z_i^+\in\R^{n\times c_i^+}$ and $Z_i^-\in\R^{n\times c_i^-}$ are the correction factors. The two equations can share computations associated with the coefficient matrix (see Section~\ref{sect:mp-ir-cost}).
In practice, it is often necessary to impose a rank truncation in the residual factorization~\cite{bequ99}, so eigenvalues of magnitude smaller than a relative threshold $0<\eta_r<1$ are discarded.
This reduces the algorithmic cost by potentially removing redundant dimensions in the correction factors.
For example, the authors of~\cite{bekq11} have observe that $\eta_r=10^{-4}$ in general works well for their algorithm.
The overall residual factorization scheme is presented as Fragment~\ref{alg.res-fac-chol}, where we use double subscripts to denote the truncated eigenvector and eigenvalue matrices.

In terms of the full solution matrices, the update of the solution factors corresponds to
\[
 X_{i+1}:=X_i+X_i^+-X_i^- -\Delta X_{i+1},
\]
where the perturbation matrix $\Delta X_{i+1}$ accounts for the projection and rank truncation of the unprojected update $X_{i+1}^{bp}:=X_i+X_i^+-X_i^-$. 
Since the sought solution $X$ is positive semidefinite, we first project $X_{i+1}^{bp}$ onto the convex cone of PSD matrices by discarding its negative eigenvalues, and then truncate small positive eigenvalues using a relative tolerance $0<\eta_s<1$.
This projection and truncation can be applied implicitly, in a similar manner to the residual factorization discussed above, without forming an $n\times n$ matrix. Writing
\[
 G_i=\begin{bmatrix}Z_i&Z_i^+&Z_i^-\end{bmatrix},\qquad
 J_i=\blkdiag(I_{c_i},I_{c_i^+},-I_{c_i^-}),
\]
we have $X_{i+1}^{bp}=G_iJ_iG_i^T$. A thin QR factorization $G_i=V_i\Gamma_i$ reduces the projection and truncation to the small kernel matrix $K_i=\Gamma_iJ_i\Gamma_i^T$, yielding $Z_{i+1}=V_i\Theta_{i,t}^+(\Sigma_{i,t}^+)^{1/2}$. 
Here $\Sigma_{i,t}^+$ and $\Theta_{i,t}^+$ contain the retained positive eigenvalues and corresponding eigenvectors of $K_i$. 
The spectral decomposition and truncation are analogous to those in Fragment~\ref{alg.res-fac-chol}; see Fragment~\ref{alg.sol-update-chol} for the complete solution update.

\begin{codefragment}[t]
 \caption{Residual factorization of Cholesky-type solution factors.}\label{alg.res-fac-chol}
 \SetAlgoVlined\nonl
 \Fn{\textsc{ResFacChol}}{
  \Parameter{Relative residual truncation tolerance $0<\eta_r<1$}
  \KwIn{$A$, $L$, $Z_i$}
  \KwOut{$L_i^+$, $L_i^-$}
  $F_i=\begin{bmatrix}Z_i&AZ_i&L\end{bmatrix}$\;
  Compute a thin QR decomposition $F_i=U_iT_i$\;
  Form $H_i=T_iP_iT_i^T$ and compute $H_i=Q_i\Lambda_iQ_i^T$\label{alg.res-fac-chol-line.spec.dec}\;
  $J^+=\{j:\lambda_j>0,\ \lambda_j\ge\eta_r\normt{\Lambda_i}\}$,
  $J^-=\{j:\lambda_j<0,\ \lambda_j\le-\eta_r\normt{\Lambda_i}\}$\;
  $\Lambda_{i,t}^{\pm}=\diag(\lambda_j)_{j\in J^{\pm}}$, $Q_{i,t}^{\pm}=Q_i(:,J^{\pm})$\;
  $L_i^+=U_iQ_{i,t}^+(\Lambda_{i,t}^+)^{1/2}$,
  $L_i^-=U_iQ_{i,t}^-(-\Lambda_{i,t}^-)^{1/2}$\;
 }
\end{codefragment}

\begin{codefragment}[t]
 \caption{Cholesky-type solution factor update.}\label{alg.sol-update-chol}
 \SetAlgoVlined\nonl
 \Fn{\textsc{SolUptChol}}{
  \Parameter{Relative solution truncation tolerance $0<\eta_s<1$}
  \KwIn{$Z_i$, $Z_i^+$, $Z_i^-$}
  \KwOut{$Z_{i+1}$}
  $G_i=\begin{bmatrix}Z_i&Z_i^+&Z_i^-\end{bmatrix}$,
  $J_i=\blkdiag(I_{c_i},I_{c_i^+},-I_{c_i^-})$\;
  Compute a thin QR decomposition $G_i=V_i\Gamma_i$\;
  Form $K_i=\Gamma_iJ_i\Gamma_i^T$ and compute $K_i=\Theta_i\Sigma_i\Theta_i^T$\label{alg.sol-fac-chol-line.spec.dec}\;
  $J^+=\{j:\sigma_j>0,\ \sigma_j\ge\eta_s\normt{\Sigma_i}\}$,
  $\Sigma_{i,t}^+=\diag(\sigma_j)_{j\in J^+}$,
  $\Theta_{i,t}^+=\Theta_i(:,J^+)$\;
  $Z_{i+1}=V_i\Theta_{i,t}^+(\Sigma_{i,t}^+)^{1/2}$\;
 }
\end{codefragment}

\begin{algorithm2e}[t]
 \caption{Mixed-precision Cholesky-type IR framework.}\label{alg.ir3-lyap-chol}
 \SetAlgoVlined
 \Parameter{Relative convergence tolerance $\tau_I>0$, maximal number of corrections $i_{\max}\in\mathbb N^+$, and precisions $u_s\ge u\ge u_r,u_c>0$}
 \KwIn{$A\in\R^{n\times n}$, $L\in\R^{n\times m}$}
 \KwOut{Solution factor $Z$ such that $X\approx ZZ^T$}
 Solve $AX+XA^T+LL^T=0$ at precision $u_s$ and store $Z_1$ at precision $u$\;
 $Z=Z_1$\;
 \For{$i\gets1$ \KwTo $i_{\max}$}{
  Evaluate $[L_i^+,L_i^-]=\textsc{ResFacChol}(A,L,Z_i)$ using Fragment~\ref{alg.res-fac-chol} at precision $u_r$\label{alg.line.ir3-lyap.resfac.chol}\;
  \If{$\operatorname{res}(Z_i)\le\tau_I$}{\textbf{break}\;}
  Round $L_i^+,L_i^-$ to precision $u_s$\;
  Solve~\eqref{eq:lyap-cor-chol} at precision $u_s$ and store $Z_i^+,Z_i^-$ at precision $u$\label{alg.line.ir3-lyap.solver.chol}\;
  Evaluate $Z_{i+1}=\textsc{SolUptChol}(Z_i,Z_i^+,Z_i^-)$ using Fragment~\ref{alg.sol-update-chol} at precision $u_c$\label{alg.line.ir3-lyap.solupt.chol}\;
  $Z=Z_{i+1}$\;
 }
\end{algorithm2e}

The mixed-precision IR framework, presented as Algorithm~\ref{alg.ir3-lyap-chol}, is essentially a \textit{projected} fixed-point iteration for solving the low-rank Lyapunov equation~\eqref{eq:lyap}.
The stopping criterion uses the relative residual $\operatorname{res}(Z_i)$ defined in~\eqref{eq:res-fac-sol} below.
The precisions are parameterized by the working precision $u$, the solver precision $u_s\ge u$, the residual factorization precision $u_r\le u$, and the solution update and projection precision $u_c\le u$. Initial and correction factors are stored at working precision $u$, while updated factors are retained at precision $u_c$ when $u_c<u$.
This basic precision setting largely aligns with the choice of precisions in the IR framework for the solution of linear systems~\cite[Alg.~1.1]{cahi18}.

\subsection{A new $\ldlt$-type IR variant}\label{sect:ldlt-ir}

In iterative solvers for low-rank Lyapunov equations, an alternative is the $\ldlt$-type formulation, which has been extensively used in ADI-based solvers~\cite{lms15}.
The solution is represented as $X=ZYZ^T$, where $Z$ is tall-and-skinny and $Y$ is symmetric and block diagonal.
This formulation accommodates any indefiniteness in the residuals and correction solutions by permitting negative eigenvalues in the diagonal blocks of $Y$, thereby allowing each correction to be computed by solving a single Lyapunov equation~\eqref{eq:lyap-cor-ldlt}.

We keep the notation consistent between the Cholesky-type and $\ldlt$-type IR schemes, to show the correspondence of the quantities and reduce repetition in the error analysis.
The residual of an approximated solution $X_i=Z_iY_iZ_i^T$, with $Z_i\in\R^{n\times c_i}$ and $Y_i\in\R^{c_i\times c_i}$, to~\eqref{eq:lyap} (with $W=LSL^T$) takes the form
\begin{equation}\label{eq:lyap-res-ldlt}
 \mathcal R(Z_i,Y_i):=AZ_iY_iZ_i^T+Z_iY_iZ_i^TA^T+LSL^T.
\end{equation}
The residual can be written in factored form as
\begin{equation}\label{eq:res-Zi-ldlt}
 \mathcal R(Z_i,Y_i)=F_iN_iF_i^T,\quad
 F_i=\begin{bmatrix}Z_i&AZ_i&L\end{bmatrix},\quad
 N_i=\begin{bmatrix}0&Y_i&0\\Y_i&0&0\\0&0&S\end{bmatrix}.
\end{equation}
A thin QR factorization $F_i=U_iT_i$ and a spectral decomposition of $H_i=T_iN_iT_i^T$ yield factors of the residual without forming an $n\times n$ matrix. 
As in the Cholesky-type residual factorization, truncating eigenvalues of small magnitude limits the rank of the residual factors and can improve efficiency. 
In the $\ldlt$-type formulation, the retained eigenvalues
are stored together in $S_i^\Delta$, rather than used to form separate factors for the PSD and NSD residual parts; cf.~\eqref{eq:res-fac}.
The overall scheme is stated as Fragment~\ref{alg.res-fac-ldlt}, which returns the $\ldlt$-type factors $L_i^{\Delta}$ and $S_i^{\Delta}$ of the truncated residual $L_i^{\Delta} S_i^{\Delta} (L_i^{\Delta})^T$.

\begin{codefragment}[t]
 \caption{Residual factorization of $\ldlt$-type solution factors.}\label{alg.res-fac-ldlt}
 \SetAlgoVlined\nonl
 \Fn{\textsc{ResFacLdlt}}{
  \Parameter{Relative residual truncation tolerance $0<\eta_r<1$}
  \KwIn{$A$, $L$, $S$, $Z_i$, $Y_i$}
  \KwOut{$L_i^\Delta$, $S_i^\Delta$}
  $F_i=\begin{bmatrix}Z_i&AZ_i&L\end{bmatrix}$\;
  $N_i=\begin{bmatrix}0&Y_i&0\\Y_i&0&0\\0&0&S\end{bmatrix}$\;
  Compute a thin QR decomposition $F_i=U_iT_i$\;
  Form $H_i=T_iN_iT_i^T$ and compute $H_i=Q_i\Lambda_iQ_i^T$\label{alg.res-fac-ldlt.line-spec-dec}\;
  $J^t=\{j:|\lambda_j|\ge\eta_r\normt{\Lambda_i}\}$,
  $S_i^\Delta=\diag(\lambda_j)_{j\in J^t}$, $Q_{i,t}=Q_i(:,J^t)$\label{alg.res-fac-ldlt.line-Si}\;
  $L_i^\Delta=U_iQ_{i,t}$\label{alg.res-fac-ldlt.line-Li}\;
 }
\end{codefragment}

The correction equation is
\begin{equation}\label{eq:lyap-cor-ldlt}
 AX_i^\Delta+X_i^\Delta A^T+L_i^\Delta S_i^\Delta(L_i^\Delta)^T=0,
 \qquad X_i^\Delta=Z_i^\Delta Y_i^\Delta(Z_i^\Delta)^T.
\end{equation}
Upon solving~\eqref{eq:lyap-cor-ldlt}, we obtain the factors $Z_i^\Delta\in\R^{n\times c_i^\Delta}$ and $Y_i^\Delta\in\R^{c_i^\Delta\times c_i^\Delta}$ of the solution increment $X_i^\Delta$.
The full solution update is
\[
 X_{i+1}:=X_i+X_i^\Delta-\Delta X_{i+1},
\]
where, as in the Cholesky-type scheme, $\Delta X_{i+1}$ accounts for PSD projection and rank truncation of $X_{i+1}^{bp}:=X_i+X_i^\Delta$. Writing
\[
 G_i=\begin{bmatrix}Z_i&Z_i^\Delta\end{bmatrix},\qquad
 \Upsilon_i=\blkdiag(Y_i,Y_i^\Delta),
\]
we have $X_{i+1}^{bp}=G_i\Upsilon_iG_i^T$. A thin QR factorization $G_i=V_i\Gamma_i$ again reduces the projection and truncation to the small kernel matrix $K_i=\Gamma_i\Upsilon_i\Gamma_i^T$, retaining positive eigenvalues at the relative threshold $\eta_s\normt{\Sigma_i}$. The updated factors are $Z_{i+1}=V_i\Theta_{i,t}^+$ and $Y_{i+1}=\Sigma_{i,t}^+$; see Fragment~\ref{alg.sol-update-ldlt}.

\begin{codefragment}[t]
 \caption{$\ldlt$-type solution factor update.}\label{alg.sol-update-ldlt}
 \SetAlgoVlined\nonl
 \Fn{\textsc{SolUptLdlt}}{
  \Parameter{Relative solution truncation tolerance $0<\eta_s<1$}
  \KwIn{$Z_i$, $Y_i$, $Z_i^\Delta$, $Y_i^\Delta$}
  \KwOut{$Z_{i+1}$, $Y_{i+1}$}
  $G_i=\begin{bmatrix}Z_i&Z_i^\Delta\end{bmatrix}$,
  $\Upsilon_i=\blkdiag(Y_i,Y_i^\Delta)$\;
  Compute a thin QR decomposition $G_i=V_i\Gamma_i$\;
  Form $K_i=\Gamma_i\Upsilon_i\Gamma_i^T$ and compute $K_i=\Theta_i\Sigma_i\Theta_i^T$\label{alg.sol-update-ldlt-line-spec.dec}\;
  $J^+=\{j:\sigma_j>0,\ \sigma_j\ge\eta_s\normt{\Sigma_i}\}$,
  $Y_{i+1}=\diag(\sigma_j)_{j\in J^+}$,
  $\Theta_{i,t}^+=\Theta_i(:,J^+)$\;
  $Z_{i+1}=V_i\Theta_{i,t}^+$\;
 }
\end{codefragment}

\begin{algorithm2e}[t]
 \caption{Mixed-precision $\ldlt$-type IR framework.}\label{alg.ir3-lyap-ldlt}
 \SetAlgoVlined
 \Parameter{Relative convergence tolerance $\tau_I>0$, maximal number of corrections $i_{\max}\in\mathbb N^+$, and precisions $u_s\ge u\ge u_r,u_c>0$}
 \KwIn{$A\in\R^{n\times n}$, $L\in\R^{n\times m}$, and symmetric PSD $S\in\R^{m\times m}$}
 \KwOut{Solution factors $Z$, $Y$ such that $X\approx ZYZ^T$}
 Solve $AX+XA^T+LSL^T=0$ at precision $u_s$ and store $Z_1,Y_1$ at precision $u$\;
 $Z=Z_1$, $Y=Y_1$\;
 \For{$i\gets1$ \KwTo $i_{\max}$}{
  Evaluate $[L_i^\Delta,S_i^\Delta]=\textsc{ResFacLdlt}(A,L,S,Z_i,Y_i)$ using Fragment~\ref{alg.res-fac-ldlt} at precision $u_r$\;
  \If{$\operatorname{res}(Z_i,Y_i)\le\tau_I$}{\textbf{break}\;}
  Round $L_i^\Delta,S_i^\Delta$ to precision $u_s$\;
  Solve~\eqref{eq:lyap-cor-ldlt} at precision $u_s$ and store $Z_i^\Delta,Y_i^\Delta$ at precision $u$\label{alg.line.ir3-lyap.solver.ldlt}\;
  Evaluate $[Z_{i+1},Y_{i+1}]=\textsc{SolUptLdlt}(Z_i,Y_i,Z_i^\Delta,Y_i^\Delta)$ using Fragment~\ref{alg.sol-update-ldlt} at precision $u_c$\;
  $Z=Z_{i+1}$, $Y=Y_{i+1}$\;
 }
\end{algorithm2e}

The mixed-precision $\ldlt$-type IR framework is presented in Algorithm~\ref{alg.ir3-lyap-ldlt}, with the same precision assignments as in the Cholesky-type scheme. 

\subsection{Rounding error analysis of $\ldlt$-type IR}\label{sect:err-anal-ir-ldlt}

\paragraph{Assumptions and notation.}
We carry out a \emph{first-order} rounding error analysis of Algorithm~\ref{alg.ir3-lyap-ldlt} in this subsection.
We use the standard model of floating-point arithmetic~\cite[sect.~2.2]{high02}
\begin{equation}\label{eq:fpmodel}
	\fl(x \hspace{2pt}\texttt{op}\hspace{2pt} y) = (x \hspace{2pt}\texttt{op}\hspace{2pt} y) (1 + \delta), \quad  |\delta| \le u,
\end{equation}
where $u$ denotes the unit roundoff associated with the floating-point number system, $x$ and $y$ are floating-point numbers, and \texttt{op} denotes
addition, subtraction, multiplication, or division.
We assume rounding to nearest in the error model and do not consider overflow and underflow. It is understood that inputs to the working precision are exactly representable.

We use the phrase ``precision $u$'' (perhaps with subscripts) to indicate a floating-point arithmetic with unit roundoff $u$.
The hats denote quantities computed in floating-point arithmetic, and $\fl_r(\cdot)$ is used to denote the computed quantity of an arithmetic process performed in precision $u_r$.
Given an integer $n$, we define 
\[
	\gamma_n=nu/(1-nu), \qquad \gammat_n= cnu / (1-cnu),
\]
where $c$ is a small integer constant whose exact value is unimportant. When $\gamma_n$ or $\gammat_n$ carries a superscript, that superscript denotes the index of the corresponding $u$ appearing as a subscript; for example $\gamma_n^s=nu_s/(1-nu_s)$. 

The accuracy is measured with respect to the full solution approximants, which are not formed by the algorithm. Therefore, the approximants to the full solution or its increments can be thought of as computed \textit{exactly} from the approximated factors, namely,
\[
 \wh X_i=\wh Z_i\wh Y_i\wh Z_i^T,\qquad
 \wh X_i^\Delta=\wh Z_i^\Delta\wh Y_i^\Delta(\wh Z_i^\Delta)^T.
\]
Similarly, since the residual~\eqref{eq:lyap-res-ldlt} is never explicitly computed, we define the computed residual exactly from the rounded factors supplied to the solver:
\[
 \wh{\mathcal R}(\wh Z_i,\wh Y_i)
 =\wh L_i^\Delta \wh S_i^\Delta \wh (L_i^\Delta)^T,
\]
where $\wh L_i^\Delta,\wh S_i^\Delta$ denote the factors computed at precision $u_r$.

For each iteration, let $F_i$ and $N_i$ be given by~\eqref{eq:res-Zi-ldlt} with $Z_i=\wh Z_i$ and $Y_i=\wh Y_i$, and set
\[
 \begin{gathered}
 F_i=U_iT_i,\qquad H_i=T_iN_iT_i^T,\\
 G_i:=\begin{bmatrix}\wh Z_i&\wh Z_i^\Delta\end{bmatrix},\qquad
 \Upsilon_i:=\blkdiag(\wh Y_i,\wh Y_i^\Delta),\\
 G_i=V_i\Gamma_i,\qquad K_i=\Gamma_i\Upsilon_i\Gamma_i^T,
 \end{gathered}
\]
where $U_i$ and $V_i$ in the  two thin QR factorizations have orthonormal columns. During the iterations, we define the dimensions and ranks
\[
 \begin{aligned}
 m_i&:=2c_i+m,\qquad &r_i&:=\rank(F_i)=\rank(T_i)\le m_i,\\
 p_i&:=c_i+c_i^\Delta, &q_i&:=\rank(G_i)=\rank(\Gamma_i)\le p_i.
 \end{aligned}
\]
The cancellation ratios used below are defined by
\begin{equation}\label{eq:assump-c1}
 \normF{|A||\wh Z_i|}=b_1\normF{F_i},
\end{equation}
where the factor $b_1\equiv b_1(n,i,A,L)$ is step dependent (as are $b_2$ and $b_3$ below) and essentially characterizes the stability of the matrix multiplication $A\wh Z_i$ with respect to potential numerical cancellation; for example, $b_1\in O(1)$ if $\normF{|A||\wh Z_i|}\approx\normF{A\wh Z_i}$.
Similarly,
\begin{equation}\label{eq:assump-c2-ldlt}
 \normF{N_i}\normF{T_iT_i^T}=b_2\normF{T_iN_iT_i^T},
\end{equation}
where the factor $b_2\equiv b_2(m,i,A,L)\ge 1$ is large if there is significant cancellation in forming the product $H_i=T_iN_iT_i^T$.
The corresponding ratio $b_3\equiv b_3(m,i,A,L)\ge 1$ in the solution update is defined by
\begin{equation}\label{eq:assump-c3-ldlt}
 \normF{\Upsilon_i}\normF{\Gamma_i\Gamma_i^T}
 =b_3\normF{\Gamma_i\Upsilon_i\Gamma_i^T}.
\end{equation}
For the analysis, we use the following assumptions.
\begin{enumerate}
 \renewcommand{\theenumi}{A\arabic{enumi}}
 \renewcommand{\labelenumi}{(\theenumi)}
 \item\label{assump:ldlt-solver}
We assume that the solver used on Line~\ref{alg.line.ir3-lyap.solver.ldlt} of Algorithm~\ref{alg.ir3-lyap-ldlt} produces computed solution factors $\wh Z_i^\Delta$ and $\wh Y_i^\Delta$ satisfying the correction equation such that
\begin{equation}\label{eq:assump-ldlt}
 \normF{\L(\wh X_i^\Delta)+\wh{\mathcal R}(\wh Z_i,\wh Y_i)}
 \le u_s\bigl(d_1\normF{\L}\normF{\wh X_i^\Delta}
              +d_2\normF{\wh{\mathcal R}(\wh Z_i,\wh Y_i)}\bigr),
\end{equation}
where the two constants $d_1,d_2>0$ depend on $A$, $\wh{\mathcal R}(\wh Z_i,\wh Y_i)$, the dimension $n$, as well as the solver precision $u_s$.
The assumption implies that the normwise relative residual of the computed solution to
$\L(X_i^\Delta)+\wh{\mathcal R}(\wh Z_i,\wh Y_i)=0$ is of order at most $\max(d_1,d_2)u_s$.
 \item\label{assump:ldlt-qr}
 The thin QR factorizations are computed by the Householder QR algorithm, and the dimensions and roundoff factors satisfy
 \[
 m_i,p_i\ll n,\qquad
 \sqrt{m_i}\gammat_{m_in}^r<1,\qquad
 \sqrt{p_i}\gammat_{p_in}^c<1.
 \]
 \item\label{assump:ldlt-truncation}
 The rank truncation tolerances satisfy $0<\eta_r,\eta_s\ll1$.
 For the negative spectral part $K_i^-$ discarded by the solution PSD projection, defined explicitly in~\eqref{eq: Ki-0-minus}, we assume
 \[
 \normF{K_i^-}=\eta_n\normF{K_i},\qquad 0\le\eta_n\ll1.
 \]
 \item\label{assump:ldlt-solver-smallness}
 The Lyapunov operator $\L(X)=AX+XA^T$ on $\R^{n\times n}$ is nonsingular, and, for the two constants $d_1,d_2>0$ from \eqref{eq:assump-ldlt},
 \[
 d_1\kappa_F(\L)u_s<1,\qquad
 \psi:=(d_1\kappa_F(\L)+d_2)u_s\ll1.
 \]
\end{enumerate}

\paragraph{Residual factorization.}
We start with analyzing the rounding errors of Fragment~\ref{alg.res-fac-ldlt}.
The computed tall-and-skinny factor $\wh F_i$ satisfies
\begin{equation}\label{eq:Fihat-def}
 \wh F_i=\begin{bmatrix}\wh Z_i&A\wh Z_i+\Delta F_i^{(2)}&L\end{bmatrix},
 \qquad |\Delta F_i^{(2)}|\le\gamma_n^r|A||\wh Z_i|.
\end{equation}
Under assumption~\ref{assump:ldlt-qr}, the thin Householder QR factorization of $\wh F_i$ satisfies~\cite[sect.~19.3, eq.~(19.13)]{high02}
\begin{equation}\label{eq:Fihat-qr-bound}
 \wh U_i\wh T_i=\wh F_i+\Delta\widetilde F_i,
 \qquad \normF{\Delta\widetilde F_i}\le\sqrt{m_i}\gammat_{m_in}^r\normF{\wh F_i},
\end{equation}
where, for some $\widetilde U_i$ with orthonormal columns,
\begin{equation}\label{eq:Uhat-def}
 \wh U_i=\widetilde U_i+\Delta U_i,
 \qquad \normF{\Delta U_i}\le\sqrt{m_i}\gammat_{m_in}^r
 \lesssim m_i^{3/2}\gammat_n^r.
\end{equation}
Here $\widetilde U_i$ need not coincide with the exact QR factor $U_i$ of $F_i$; thus this backward relation does not imply a small forward error in the individual QR factors.

Writing $\wh F_i+\Delta\widetilde F_i=:F_i+\Delta F_i$ and combining~\eqref{eq:Fihat-def}--\eqref{eq:Fihat-qr-bound} and~\eqref{eq:assump-c1}, we obtain the first-order rounding error bound
\begin{equation}\label{eq:errbnd-deltaT}
 \wh U_i\wh T_i=F_i+\Delta F_i,
 \qquad \normF{\Delta F_i}\lesssim
 (b_1+m_i^{3/2})\gammat_n^r\normF{F_i}.
\end{equation}
Since $\widetilde U_i$ has orthonormal columns, the triangle and reverse triangle inequalities give
\begin{equation}\label{eq:Uhat-norm-bounds}
 \normt{\wh U_i}\le1+\normt{\Delta U_i},\qquad
 \sigma_{\min}(\wh U_i)\ge1-\normt{\Delta U_i}>0,
\end{equation}
and consequently, one can show that
\begin{equation}\label{eq:That-norm-bound}
 \normF{\wh T_i}
 \le\frac{\normF{\wh U_i\wh T_i}}{1-\normt{\Delta U_i}}
 =\frac{\normF{F_i+\Delta F_i}}{1-\normt{\Delta U_i}}.
\end{equation}

For Line~\ref{alg.res-fac-ldlt.line-spec-dec} of Fragment~\ref{alg.res-fac-ldlt}, define $\wh H_i:=\fl_r((\wh T_iN_i)\wh T_i^T)$, where $N_i$ is assembled from the stored small factors. Therefore, we have, to first order,
\begin{equation}\label{eq:Hhat-formation-ldlt}
 \wh H_i=\wh T_iN_i\wh T_i^T+\Delta\widetilde H_i,
 \qquad |\Delta\widetilde H_i|\lesssim
 2\gammat_{m_i}^r|\wh T_i||N_i||\wh T_i|^T,
\end{equation}
where the factor $2$ accounts for the two matrix multiplications performed sequentially, and then
\begin{equation}\label{eq:deltaH-def-ldlt}
 \Delta H_i:=\wh H_i-H_i
 =\wh T_iN_i\wh T_i^T-T_iN_iT_i^T+\Delta\widetilde H_i.
\end{equation}
Substituting~\eqref{eq:errbnd-deltaT} into~\eqref{eq:Hhat-formation-ldlt} yields the following perturbation of the \emph{reconstructed residual}:
\[
 \wh U_i\wh H_i\wh U_i^T-F_iN_iF_i^T
 =F_iN_i\Delta F_i^T+\Delta F_iN_iF_i^T
   +\Delta F_iN_i\Delta F_i^T
 +\wh U_i\Delta\widetilde H_i\wh U_i^T.
\]
Using~\eqref{eq:errbnd-deltaT}--\eqref{eq:Hhat-formation-ldlt}, one can get the first-order rounding error bound
\begin{equation}\label{eq:res-perturb-bound}
 \normF{\wh U_i\wh H_i\wh U_i^T-F_iN_iF_i^T}
 \lesssim
 \bigl(2(b_1+m_i^{3/2})\gammat_n^r+2\gammat_{m_i}^r\bigr)
 \normF{N_i}\normF{F_i}^2.
\end{equation}
Applying the Cauchy--Schwarz inequality to the squared singular values of $T_i$ gives
\[
 \normF{T_iT_i^T}\le\normF{T_i}^2
 \le\sqrt{r_i}\normF{T_iT_i^T}.
\]
Since $U_i^TU_i=I$, we have $\normF{F_i}=\normF{T_i}$ and $\normF{F_iN_iF_i^T}=\normF{H_i}$.
Using~\eqref{eq:assump-c2-ldlt}, the perturbation bound~\eqref{eq:res-perturb-bound} becomes
\begin{equation}\label{eq:reconstructed-H-bnd}
 \normF{\wh U_i\wh H_i\wh U_i^T-F_iN_iF_i^T}
 \lesssim b_2\sqrt{r_i}
 \bigl(2(b_1+m_i^{3/2})\gammat_n^r+2\gammat_{m_i}^r\bigr)\normF{H_i}.
\end{equation}
Using the singular-value bound~\eqref{eq:Uhat-norm-bounds}, the triangle inequality, and the error bounds~\eqref{eq:Uhat-def} and~\eqref{eq:reconstructed-H-bnd}, we obtain
\begin{equation}\label{eq:Hhat-norm-bound}
 \begin{aligned}
 \normF{\wh H_i}
&
\le\frac{\normF{H_i}+\normF{\wh U_i\wh H_i\wh U_i^T-F_iN_iF_i^T}}
 {(1-\normt{\Delta U_i})^2}\\
 &\lesssim\Bigl[1+
 2\bigl(b_2\sqrt{r_i}(b_1+m_i^{3/2})+m_i^{3/2}\bigr)\gammat_n^r
 +2b_2\sqrt{r_i}\gammat_{m_i}^r\Bigr]\normF{H_i},
 \end{aligned}
\end{equation}
Let $\lambda_j$ denote an eigenvalue of $H_i$, and let $\wt\lambda_j$ denote an exact eigenvalue of $\wh H_i$.
Applying Weyl's inequality to the reconstructed residuals, with their eigenvalues ordered consistently, gives,
\begin{equation}\label{eq:weyl-ineq-conseq}
 \begin{aligned}
 |\wt\lambda_j-\lambda_j|
 &\lesssim\normF{\widetilde U_i\wh H_i\widetilde U_i^T-U_iH_iU_i^T}\\
 &\lesssim\normF{\wh U_i\wh H_i\wh U_i^T-F_iN_iF_i^T}
+m_i^{3/2}\gammat_n^r\normF{\wh H_i}\\
 &\lesssim\Bigl[b_2\sqrt{r_i}
  \bigl(2(b_1+m_i^{3/2})\gammat_n^r+2\gammat_{m_i}^r\bigr)
 +m_i^{3/2}\gammat_n^r\Bigr]\normF{H_i},
\end{aligned}
\end{equation}
where the second inequality uses~\eqref{eq:Uhat-def} and the last uses~\eqref{eq:reconstructed-H-bnd} and~\eqref{eq:Hhat-norm-bound}.

The spectral decomposition on Line~\ref{alg.res-fac-ldlt.line-spec-dec} of Fragment~\ref{alg.res-fac-ldlt} is usually computed by the symmetric QR algorithm or the divide-and-conquer algorithm after a tridiagonalization~\cite[sect.~5.3]{demm97},~\cite[sect.~8.3]{gova13},~\cite{guei95}, and the computed eigensystem of $\wh H_i$ satisfies
\begin{equation}\label{eq:Hhat-eig-bnd}
 \wh Q_i\wh\Lambda_i\wh Q_i^T=\wh H_i+\Delta\wh H_i,
 \qquad \normF{\Delta\wh H_i}\lesssim\gammat_{m_i}^r\normF{\wh H_i},
\end{equation}
where $\wh\Lambda_i=\diag(\wh\lambda_j)$ contains the \textit{computed} eigenvalues of $\wh H_i$ and $\wh Q_i$ is (numerically) orthogonal in precision $u_r$.
The expression~\eqref{eq:Hhat-eig-bnd} means that the eigensystem of $\wh H_i$ is computed backward stably.
It follows from this eigensolver bound and~\eqref{eq:Hhat-norm-bound}--\eqref{eq:weyl-ineq-conseq} that, for any computed eigenvalue $\wh\lambda_j$ of $\wh H_i$, with the same eigenvalue ordering,
\begin{equation}\label{eq:H-respec-eig-bnd-ldlt}
 \begin{aligned}
 |\wh\lambda_j-\lambda_j|
 &\le|\wh\lambda_j-\wt\lambda_j|+|\wt\lambda_j-\lambda_j|
 \lesssim\gammat_{m_i}^r\normF{H_i+\Delta H_i}
       +|\wt\lambda_j-\lambda_j|\\
 &\lesssim\Bigl[b_2\sqrt{r_i}
  \bigl(2(b_1+m_i^{3/2})\gammat_n^r+2\gammat_{m_i}^r\bigr)
 +m_i^{3/2}\gammat_n^r+\gammat_{m_i}^r\Bigr]\normF{H_i},
 \end{aligned}
\end{equation}
where $\normF{H_i}=(\sum_j|\lambda_j|^2)^{1/2}\le\sum_j|\lambda_j|$.
The bound~\eqref{eq:H-respec-eig-bnd-ldlt} implies that the eigenvalues of large magnitude of the residual $\mathcal R(\wh Z_i,\wh Y_i)$ will be computed to high relative accuracy, if the rounding error factors are small (necessarily when the factors $b_1$ and $b_2$ are of moderate size).
But, in any case, there is no guarantee for the relative accuracy of computed eigenvalues of small magnitude.

Recall that, when evaluated in exact arithmetic, the residual of the computed factors $(\wh Z_i,\wh Y_i)$ satisfies
\begin{equation}\label{eq:res-Zi-def}
 \mathcal R(\wh Z_i,\wh Y_i)= L_i^\Delta  S_i^\Delta (L_i^\Delta)^T=F_iN_iF_i^T=U_iH_iU_i^T,
\end{equation}
where the residual is not explicitly formed in Algorithm~\ref{alg.ir3-lyap-ldlt}.
Therefore, to quantify the rounding errors in the computed residual $\wh{\mathcal R}(\wh Z_i,\wh Y_i)$, we only need to consider the inexact computation of the factors $L_i^\Delta,S_i^\Delta$ and its effect on the residual, and so we can write the computed residual as
\[
 \wh{\mathcal R}(\wh Z_i,\wh Y_i)
 =\wh L_i^\Delta\wh S_i^\Delta(\wh L_i^\Delta)^T+\Delta R_i^s,
\]
where $\Delta R_i^s$ accounts for the error from rounding the factors $\wh L_i^\Delta$ and $\wh S_i^\Delta$ down to precision $u_s$.
The precision conversion incurs relative componentwise perturbations to these factors bounded above by $u_s$.
Since $\wh S_i^\Delta$ is diagonal and $\wh L_i^\Delta$ has numerically orthonormal columns,
\[
 \normF{\bigl(\fl_s(\wh L_i^\Delta)-\wh L_i^\Delta\bigr)\wh S_i^\Delta}
 \le u_s\normF{|\wh L_i^\Delta||\wh S_i^\Delta|}
 \lesssim u_s\normF{\wh S_i^\Delta}.
\]
The perturbations therefore satisfy
\begin{equation}\label{eq:res-conversion-bound}
 \normF{\Delta R_i^s}\lesssim
 3u_s\normF{\wh{\mathcal R}(\wh Z_i,\wh Y_i)}.
\end{equation}

On completion of Fragment~\ref{alg.res-fac-ldlt}, we have
\[
 \begin{gathered}
 \wh S_i^\Delta=\wh\Lambda_{i,t},\qquad
 \wh L_i^\Delta=\fl_r(\wh U_i\wh Q_{i,t}),\\
 \wh Q_i=\begin{bmatrix}\wh Q_{i,t}&\wh Q_{i,0}\end{bmatrix},\qquad
 \wh\Lambda_i=\blkdiag(\wh\Lambda_{i,t},\wh\Lambda_{i,0}),
 \end{gathered}
\]
where the subscripts $t$ and $0$ label the retained and discarded eigenpairs, respectively. 
It is straightforward to show that
\[
 \wh L_i^\Delta=\wh U_i\wh Q_{i,t}+\Delta\wh L_i^\Delta,
 \qquad
 |\Delta\wh L_i^\Delta|\le\gamma_{m_i}^r|\wh U_i||\wh Q_{i,t}|,
\]
and then
\begin{equation}\label{eq:residual-factor-product-ldlt}
 \wh L_i^\Delta\wh S_i^\Delta(\wh L_i^\Delta)^T
 =\wh U_i\wh Q_{i,t}\wh\Lambda_{i,t}\wh Q_{i,t}^T\wh U_i^T+\Delta R_i, \qquad
 |\Delta R_i|\lesssim2\gamma_{m_i}^r|\wh U_i||\wh Q_{i,t}||\wh\Lambda_{i,t}||\wh Q_{i,t}^T||\wh U_i^T|.
\end{equation}
A normwise bound for $\Delta R_i$ follows immediately, but a sharper bound can be obtained by exploiting numerical orthogonality of $\wh U_i$ and $\wh Q_{i,t}$, as we now show.
Since $\wh\Lambda_{i,t}$ is diagonal, we have
\begin{equation}\label{eq:residual-factor-weighted-bound-ldlt}
\begin{aligned}
	\normF{\Delta\wh L_i^\Delta\wh\Lambda_{i,t}}
	&\le\gamma_{m_i}^r\normt{|\wh U_i|}
	\normF{|\wh Q_{i,t}||\wh\Lambda_{i,t}|}
	=\gamma_{m_i}^r\normt{|\wh U_i|}
	\normF{\wh Q_{i,t}\wh\Lambda_{i,t}}\\
	&\lesssim\sqrt{m_i}\gamma_{m_i}^r\normF{\wh\Lambda_{i,t}},
\end{aligned}
\end{equation}
where we have used the inequality $\normF{AB}\le\normt{A}\normF{B}$.
Expanding the perturbed factor product in~\eqref{eq:residual-factor-product-ldlt} gives
\[
  \Delta R_i
  =\Delta\wh L_i^\Delta\wh\Lambda_{i,t}\wh Q_{i,t}^T\wh U_i^T
  +\wh U_i\wh Q_{i,t}\wh\Lambda_{i,t}
       (\Delta\wh L_i^\Delta)^T
  + \Delta\wh L_i^\Delta\wh\Lambda_{i,t}
       (\Delta\wh L_i^\Delta)^T.
\]
Omitting the quadratic term and applying~\eqref{eq:Hhat-eig-bnd} yields,
\[
 \begin{aligned}
  \normF{\Delta R_i}
  &\lesssim2\normF{\Delta\wh L_i^\Delta\wh\Lambda_{i,t}}
       \normt{\wh U_i\wh Q_{i,t}}
  \lesssim2\sqrt{m_i}\gamma_{m_i}^r\normF{\wh\Lambda_{i,t}}\\
  &\le2\sqrt{m_i}\gamma_{m_i}^r\normF{\wh\Lambda_i}
   \lesssim2\sqrt{m_i}\gamma_{m_i}^r\normF{\wh H_i}.
 \end{aligned}
\]
Hence, by~\eqref{eq:Hhat-eig-bnd}, the computed residual can be rewritten as
\begin{equation*}
 \wh{\mathcal R}(\wh Z_i,\wh Y_i)
 =\wh U_i(\wh H_i+\Delta\wh H_i)\wh U_i^T
 -\wh U_i\wh Q_{i,0}\wh\Lambda_{i,0}\wh Q_{i,0}^T\wh U_i^T
 +\Delta R_i+\Delta R_i^s,\qquad 
 \normt{\wh\Lambda_{i,0}}\le\eta_r\normt{\wh\Lambda_i}.
\end{equation*}
Define $\wh H_i^0:=\wh Q_{i,0}\wh\Lambda_{i,0}\wh Q_{i,0}^T$.
As $\wh Q_{i,0}$ has numerically orthonormal columns, we have
\begin{equation}\label{eq:Hhat-truncation-bounds}
 \normt{\wh H_i^0}\lesssim2\eta_r\normt{\wh H_i},\qquad
 \normF{\wh H_i^0}\lesssim2\sqrt{m_i}\eta_r\normF{\wh H_i}.
\end{equation}
Hence, using~\eqref{eq:res-Zi-def}, we obtain the perturbation identity
\begin{align}\label{eq:res-delta-comput-Zi-def-ldlt}
 \Delta\mathcal R(\wh Z_i,\wh Y_i)
 &:=\wh{\mathcal R}(\wh Z_i,\wh Y_i)-\mathcal R(\wh Z_i,\wh Y_i) \nonumber\\
 &=\wh U_i\wh H_i\wh U_i^T-F_iN_iF_i^T
 +\wh U_i(\Delta\wh H_i-\wh H_i^0)\wh U_i^T
       +\Delta R_i+\Delta R_i^s.
\end{align}
Under the truncation assumption~\ref{assump:ldlt-truncation}, a first-order normwise bound follows immediately by using~\eqref{eq:reconstructed-H-bnd},~\eqref{eq:Hhat-norm-bound},~\eqref{eq:Hhat-eig-bnd},~\eqref{eq:res-conversion-bound}, and~\eqref{eq:Hhat-truncation-bounds}, and we have  
 \begin{align}\label{eq:res-comput-Zi-bnd-ldlt}
  \normF{\Delta\mathcal R(\wh Z_i,\wh Y_i)}
  &\lesssim\Bigl[2\sqrt{m_i}\eta_r+3u_s
  +2 b_2\sqrt{r_i}
       (b_1+m_i^{3/2})\gammat_n^r+2(b_2\sqrt{r_i} + \sqrt{m_i}) \gammat_{m_i}^r \Bigr]
       \normF{\mathcal R(\wh Z_i,\wh Y_i)} \nonumber\\
   &\approx \Bigl[
   2\sqrt{m_i}\eta_r+3u_s +2 b_2\sqrt{r_i}   (b_1+m_i^{3/2})\gammat_n^r 
   \Bigr]
   \normF{\mathcal R(\wh Z_i,\wh Y_i)}.
 \end{align}
The approximation in the second line of~\eqref{eq:res-comput-Zi-bnd-ldlt} follows from the low-rank condition in assumption~\ref{assump:ldlt-qr}, for which $\gammat_{m_i}^c\approx(m_i/n)\gammat_n^c$.
It clearly shows that the inaccuracy in the computed residual $\Delta\mathcal R(\wh Z_i,\wh Y_i)$ will be dominated by the rank truncation error when
\[
 \eta_r\gtrsim \frac{3u_s}{2\sqrt{m_i}}
 +b_2\sqrt{\frac{r_i}{m_i}}(b_1+m_i^{3/2})\gammat_n^r.
\]
\paragraph{Solution update.}
The analysis of Fragment~\ref{alg.sol-update-ldlt} is similar to the discussion above. Under assumption~\ref{assump:ldlt-qr}, the thin Householder QR factorization of $G_i$ satisfies~\cite[sect.~19.3]{high02}
\begin{equation}\label{eq:qr-Gi-hat}
 \wh V_i\wh\Gamma_i=G_i+\Delta G_i,\qquad
 \normF{\Delta G_i}\le\sqrt{p_i}\gammat_{p_in}^c\normF{G_i},
\end{equation}
and
\begin{equation}\label{eq:Vhat-def}
 \wh V_i=\widetilde V_i+\Delta V_i,\qquad
 \normF{\Delta V_i}\le\sqrt{p_i}\gammat_{p_in}^c
 \lesssim p_i^{3/2}\gammat_n^c,
\end{equation}
with $\widetilde V_i$ having orthonormal columns.
Define $\wh K_i:=\fl_c((\wh\Gamma_i\Upsilon_i)\wh\Gamma_i^T)$.
One can then show, similarly to the derivation of~\eqref{eq:deltaH-def-ldlt}--\eqref{eq:reconstructed-H-bnd} and using~\eqref{eq:assump-c3-ldlt},
\begin{equation}\label{eq:deltaK-def-ldlt}
 \begin{aligned}
  \Delta K_i&:=\wh K_i-K_i,\\
  \normF{\wh V_i\wh K_i\wh V_i^T-G_i\Upsilon_iG_i^T}
  &\lesssim b_3\sqrt{q_i}
    \bigl(2p_i^{3/2}\gammat_n^c+2\gammat_{p_i}^c\bigr)\normF{K_i},
 \end{aligned}
\end{equation}
For the spectral decomposition at Line~\ref{alg.sol-update-ldlt-line-spec.dec} of Fragment~\ref{alg.sol-update-ldlt}, the computed eigensystem of $\wh K_i$ satisfies~\cite[sect.~5.3]{demm97},~\cite[sect.~8.3]{gova13},~\cite{guei95}
\begin{equation}\label{eq:Khat-eig-bnd}
 \wh\Theta_i\wh\Sigma_i\wh\Theta_i^T=\wh K_i+\Delta\wh K_i,\qquad
 \normF{\Delta\wh K_i}\lesssim\gammat_{p_i}^c\normF{\wh K_i},
\end{equation}
where $\wh\Sigma_i$ contains the computed eigenvalues of $\wh K_i$ and $\wh\Theta_i$ is (numerically) orthogonal in precision $u_c$.
Overall, we have, to first order,
\[
 \wh Y_{i+1}=\wh\Sigma_{i,t}^+,\qquad
 \wh Z_{i+1}=\wh V_i\wh\Theta_{i,t}^++\Delta\wh Z_{i+1},\qquad
 |\Delta\wh Z_{i+1}|\le\gamma_{p_i}^c|\wh V_i||\wh\Theta_{i,t}^+|,
\]
and hence
\begin{equation}\label{eq:Xhat-i+1-ldlt}
 \begin{aligned}
  \wh X_{i+1}
  &=\wh V_i\wh\Theta_{i,t}^+\wh\Sigma_{i,t}^+
       (\wh\Theta_{i,t}^+)^T\wh V_i^T+\Delta\wh X_{i+1},\\
  \normF{\Delta\wh X_{i+1}}
  &\lesssim2\sqrt{p_i}\gamma_{p_i}^c\normF{\wh X_{i+1}},
 \end{aligned}
\end{equation}
where $\wh\Theta_{i,t}^+$ and $\wh\Sigma_{i,t}^+$ collect the eigenvectors and eigenvalues, respectively, that are retained after the rank truncation of the computed full eigensystem
\begin{equation}\label{eq:comput-full-eigsystem}
  \wh\Theta_i\wh\Sigma_i\wh\Theta_i^T
  =\begin{bmatrix}\wh\Theta_{i,t}^+&\wh\Theta_{i,0}^+&\wh\Theta_i^-
  \end{bmatrix}
  \cdot  \blkdiag(\wh\Sigma_{i,t}^+,\wh\Sigma_{i,0}^+,\wh\Sigma_i^-)
  \cdot\begin{bmatrix}\wh\Theta_{i,t}^+&\wh\Theta_{i,0}^+&\wh\Theta_i^-\end{bmatrix}^T.
\end{equation}
The bound in~\eqref{eq:Xhat-i+1-ldlt} is obtained by weighting each error column with the diagonal matrix $\wh\Sigma_{i,t}^+$, yielding the intermediate estimate
\[
 \normF{\Delta\wh Z_{i+1}\wh\Sigma_{i,t}^+}
 \le\gamma_{p_i}^c\normF{\wh V_i}
       \normF{\wh\Theta_{i,t}^+\wh\Sigma_{i,t}^+}
 \lesssim\sqrt{p_i}\gamma_{p_i}^c\normF{\wh\Sigma_{i,t}^+}.
\]
Define
\begin{equation}\label{eq: Ki-0-minus}
 K_i^0:=\wh\Theta_{i,0}^+\wh\Sigma_{i,0}^+(\wh\Theta_{i,0}^+)^T,\qquad
 K_i^-:=\wh\Theta_i^-\wh\Sigma_i^-(\wh\Theta_i^-)^T.
\end{equation}
As in~\eqref{eq:Hhat-truncation-bounds}, numerical orthonormality of $\wh\Theta_{i,0}^+$ and the eigenvalue truncation criterion give
\[
	\normt{K_i^0}\lesssim2\eta_s\normt{\wh K_i+\Delta\wh K_i},
	\qquad 
	\normF{K_i^0}\lesssim2\sqrt{p_i}\eta_s\normF{\wh K_i+\Delta\wh K_i}.
\]
The matrix $K_i^-$ corresponds to the negative semidefinite part discarded by the PSD projection; its norm satisfies assumption~\ref{assump:ldlt-truncation}.

Note that $\wh X_i+\wh X_i^\Delta=G_i\Upsilon_iG_i^T=V_iK_iV_i^T$, where the summation can be considered exact as it is performed implicitly via Fragment~\ref{alg.sol-update-ldlt}.
We have, from~\eqref{eq:Khat-eig-bnd}--\eqref{eq: Ki-0-minus},
\[
 \wh X_{i+1}=\wh V_i(\wh K_i+\Delta\wh K_i)\wh V_i^T
       -\wh V_i(K_i^0+K_i^-)\wh V_i^T+\Delta\wh X_{i+1},
\]
and then the identity of first-order perturbation
\begin{align}\label{eq:sol-upt-delta-comput-def-ldlt}
 \Delta\Xi_i&:=\wh X_{i+1}-(\wh X_i+\wh X_i^\Delta)\nonumber\\
 &=\wh V_i\wh K_i\wh V_i^T-G_i\Upsilon_iG_i^T +\wh V_i(\Delta\wh K_i-K_i^0-K_i^-)\wh V_i^T
       +\Delta\wh X_{i+1}.
\end{align}
By using~\eqref{eq:Vhat-def}--\eqref{eq:deltaK-def-ldlt},~\eqref{eq:Khat-eig-bnd}--\eqref{eq:Xhat-i+1-ldlt}, the preceding truncation bounds, and assumption~\ref{assump:ldlt-truncation}, a first-order normwise bound follows immediately:
\begin{align}\label{eq:sol-upt-delta-compute-normbnd}
 \normF{\Delta\Xi_i}
 &\lesssim\normF{\wh V_i\wh K_i\wh V_i^T-G_i\Upsilon_iG_i^T}
       +\normF{\Delta\wh K_i-K_i^0-K_i^-}
       +2\sqrt{p_i}\gamma_{p_i}^c\normF{\wh X_{i+1}}\nonumber\\
 &\lesssim\Bigl[2\sqrt{p_i}\eta_s+\eta_n
       +2b_3\sqrt{q_i}\,p_i^{3/2}\gammat_n^c
       +(2b_3\sqrt{q_i}+1+2\sqrt{p_i})\gammat_{p_i}^c\Bigr]
       \normF{\wh X_{i+1}}\nonumber\\
 &\approx\Bigl[2\sqrt{p_i}\eta_s+\eta_n
       +2b_3\sqrt{q_i}\,p_i^{3/2}\gammat_n^c\Bigr]
       \normF{\wh X_{i+1}}.
\end{align}
where the last approximation follows from the low-rank condition in assumption~\ref{assump:ldlt-qr}.

\paragraph{Residual reduction.}
Using the update identity~\eqref{eq:sol-upt-delta-comput-def-ldlt} and applying the residual definition~\eqref{eq:lyap-res-ldlt} at iterations $i+1$ and $i$, we obtain, by linearity of $\L$,
\begin{equation}\label{eq:res-Zi+1}
 \begin{aligned}
  \mathcal R(\wh Z_{i+1},\wh Y_{i+1})
  &=\L(\wh X_{i+1})-\L(X)\\
  &=\mathcal R(\wh Z_i,\wh Y_i)+\L(\wh X_i^\Delta)+\L(\Delta\Xi_i).
 \end{aligned}
\end{equation}
Defining $W_i:=\L(\wh X_i^\Delta)+\wh{\mathcal R}(\wh Z_i,\wh Y_i)$, by assumption~\eqref{eq:assump-ldlt} we have
\begin{align*}
 \normF{W_i}
 &\le u_s\Bigl(d_1\normF{\L}\normF{\wh X_i^\Delta}
       +d_2\normF{\wh{\mathcal R}(\wh Z_i,\wh Y_i)}\Bigr)\\
 &\le u_s\Bigl(d_1\kappa_F(\L)
       \bigl(\normF{\wh{\mathcal R}(\wh Z_i,\wh Y_i)}+\normF{W_i}\bigr)
       +d_2\normF{\wh{\mathcal R}(\wh Z_i,\wh Y_i)}\Bigr).
\end{align*}
Rearranging the preceding inequality gives the bound
\begin{equation}\label{eq:Wi-bnd}
 \normF{W_i}\le
 u_s\frac{d_1\kappa_F(\L)+d_2}{1-d_1\kappa_F(\L)u_s}
       \normF{\wh{\mathcal R}(\wh Z_i,\wh Y_i)}
   \approx 
   \frac{\psi}{1-d_1\kappa_F(\L)u_s}
   \normF{\mathcal R(\wh Z_i,\wh Y_i)},
\end{equation}
whose denominator is positive by assumption~\ref{assump:ldlt-solver-smallness}.
Note that the approximation in~\eqref{eq:Wi-bnd} is to first-order, using~\eqref{eq:res-delta-comput-Zi-def-ldlt}--\eqref{eq:res-comput-Zi-bnd-ldlt}.
Substituting $W_i$ back into~\eqref{eq:res-Zi+1} and using~\eqref{eq:res-delta-comput-Zi-def-ldlt} gives
\[
 \mathcal R(\wh Z_{i+1},\wh Y_{i+1})
 =W_i-\Delta\mathcal R(\wh Z_i,\wh Y_i)+\L(\Delta\Xi_i).
\]
Taking the norm on both sides and applying the bounds~\eqref{eq:res-comput-Zi-bnd-ldlt},~\eqref{eq:sol-upt-delta-compute-normbnd}, and~\eqref{eq:Wi-bnd}, we have, to first order,
\begin{align}\label{eq:common-residual-recurrence}
 \normF{\mathcal R(\wh Z_{i+1},\wh Y_{i+1})}
 &\le\normF{W_i}+\normF{\Delta\mathcal R(\wh Z_i,\wh Y_i)}
       +\normF{\L}\normF{\Delta\Xi_i}\nonumber\\
 &\lesssim \Biggl[
 \frac{\psi}{1-d_1\kappa_F(\L)u_s}
       +2\sqrt{m_i}\eta_r+3u_s
 +2b_2\sqrt{r_i}(b_1+m_i^{3/2})\gammat_n^r
       \Biggr]\normF{\mathcal R(\wh Z_i,\wh Y_i)}\nonumber\\
 &\quad{}+\Bigl[2\sqrt{p_i}\eta_s+\eta_n
       +2b_3\sqrt{q_i}\,p_i^{3/2}\gammat_n^c\Bigr]
       \normF{\L}\normF{\wh X_{i+1}}.
\end{align}

The next theorem summarizes our analysis of the behavior of the residual of Algorithm~\ref{alg.ir3-lyap-ldlt}.

\begin{theorem}\label{thm:res-converge-ldlt}
Let Algorithm~\ref{alg.ir3-lyap-ldlt} be applied to a Lyapunov equation $\L(X)+LSL^T=0$ under assumptions~\ref{assump:ldlt-solver}--\ref{assump:ldlt-solver-smallness}.
Within the first-order analysis and low-rank approximation used in~\eqref{eq:common-residual-recurrence}, if $\psi:=(d_1\kappa_F(\L)+d_2)u_s$ is sufficiently smaller than $1$ and
\begin{equation}\label{eq:phi-ldlt}
 \begin{aligned}
 \phi_i&:=\psi+2\sqrt{m_i}\eta_r+3u_s
       +2b_2\sqrt{r_i}(b_1+m_i^{3/2})\gammat_n^r,\\
 \phi&:=\sup_{i\ge1}\phi_i<1,
 \end{aligned}
\end{equation}
then the normwise residual is reduced on the $i$th iteration by a factor of at most approximately $\phi_i$.
If the supremum in the bound below is finite, the limiting residual satisfies
\[
  \limsup_{i\to\infty}\normF{\mathcal R(\wh Z_i,\wh Y_i)}
  \lesssim\frac{\normF{\L}}{1-\phi}\sup_{i\ge1}\Biggl\{
       \Bigl[2\sqrt{p_i}\eta_s+\eta_n
       +2b_3\sqrt{q_i}\,p_i^{3/2}\gammat_n^c\Bigr]
       \normF{\wh X_{i+1}}\Biggr\},
\]
where $b_1$, $b_2$, and $b_3$ are the step-dependent factors defined in~\eqref{eq:assump-c1},~\eqref{eq:assump-c2-ldlt}, and~\eqref{eq:assump-c3-ldlt}, respectively.
\end{theorem}

The residual reduction rate $\phi$ depends not only on $\psi$, which essentially concerns $\kappa_F(\L)u_s$, but also on the residual truncation parameter $\eta_r$ and the precision $u_r$ at which the residual factorization is performed.
Given that $u_r<u_s<\psi$, possibly by a large margin, the potential instability in the residual factorization at precision $u_r$, as indicated by the factors $b_1$ and $b_2$, is unlikely to affect the residual reduction rate of the algorithm.
Also, the optimal value of $\eta_r$ should satisfy $\eta_r \in O(\kappa_F(\L)u_s)$ to achieve best efficiency and balanced terms in the residual reduction rate.

Aiming at attaining a relative residual of order $nu$, one can set $\eta_s \in O(nu)$, and choose $u_c =u$ if the solution update is stable and the corresponding factor $b_3$ is of moderate size. Otherwise, higher precision should be used for the solution update, with $u_c \in O(u/b_3)$.

Table~\ref{table:ir3} lists precision combinations for Algorithm~\ref{alg.ir3-lyap-ldlt} and Algorithm~\ref{alg.ir3-lyap-chol}, together with indicative conditioning scales and target relative residuals.
These targets are conditional on assumptions~\ref{assump:ldlt-solver}--\ref{assump:ldlt-solver-smallness} and the convergence conditions of Theorem~\ref{thm:res-converge-ldlt}, with the error bound on solution update small enough to attain the stated accuracy.
For the Cholesky-type scheme, the solver assumption~\eqref{eq:assump-chol} replaces~\eqref{eq:assump-ldlt}, with the modifications described in Section~\ref{sect:chol-specialization}.

\begin{table}[t]
	\caption{Precision combinations for mixed-precision IR, conditioning bounds (two significant figures), and target relative residuals of order $nu$. Target residuals are conditional on assumptions~\ref{assump:ldlt-solver}--\ref{assump:ldlt-solver-smallness}, $\phi<1$ in~\eqref{eq:phi-ldlt}, and the bound~\eqref{eq:sol-upt-delta-compute-normbnd} for solution update being small enough for Theorem~\ref{thm:res-converge-ldlt} to yield the stated accuracy.}\label{table:ir3}
	\centering\setlength\tabcolsep{8pt}
	\begin{tabular}{c c c c}
		\toprule
		$u_s$ & $u_r=u_c=u$ & Conditioning bound $1/u_s$ & Target residual\\
		\midrule
		bf16 & \multirow{3}{*}{fp32} & $2.6\times10^2\phantom{1}$ & \multirow{3}{*}{fp32}\\
		fp16 & & $2.0\times10^3\phantom{1}$ & \\
		fp32 & & $1.7\times10^7\phantom{1}$ & \\
		\midrule
		bf16 & \multirow{4}{*}{fp64} & $2.6\times10^2\phantom{1}$ & \multirow{4}{*}{fp64}\\
		fp16 & & $2.0\times10^3\phantom{1}$ & \\
		fp32 & & $1.7\times10^7\phantom{1}$ & \\
		fp64 & & $9.0\times10^{15}$ & \\
		\bottomrule
	\end{tabular}
\end{table}

\subsection{Specialization of the analysis to Cholesky-type IR}\label{sect:chol-specialization}
The analysis in Section~\ref{sect:err-anal-ir-ldlt} applies to Algorithm~\ref{alg.ir3-lyap-chol} with
\begin{equation}\label{eq:chol-specialization-notation}
 \begin{gathered}
 N_i=P_i,\qquad \Upsilon_i=J_i,\qquad
 G_i=\begin{bmatrix}\wh Z_i&\wh Z_i^+&\wh Z_i^-\end{bmatrix},\\
 p_i=c_i+c_i^++c_i^-,\qquad r_i=\rank(F_i),\qquad q_i=\rank(G_i).
 \end{gathered}
\end{equation}
Here $\wh X_i=\wh Z_i\wh Z_i^T$, $\wh X_i^\pm=\wh Z_i^\pm(\wh Z_i^\pm)^T$, and $\wh X_i^\Delta=\wh X_i^+-\wh X_i^-$ can be viewed as exactly reconstructed from the stored factors.
Assumptions~\ref{assump:ldlt-qr}--\ref{assump:ldlt-solver-smallness} are retained with the notation in~\eqref{eq:chol-specialization-notation}, including the bound on $K_i^-$ in assumption~\ref{assump:ldlt-truncation}.
The solver assumption~\ref{assump:ldlt-solver} is replaced by the following bound for the difference of the two corrections computed on Line~\ref{alg.line.ir3-lyap.solver.chol} of Algorithm~\ref{alg.ir3-lyap-chol}:
\begin{equation}\label{eq:assump-chol}
 \begin{aligned}
 &\normF{\L(\wh X_i^+-\wh X_i^-)
   +\wh{\mathcal R}^+(\wh Z_i)-\wh{\mathcal R}^-(\wh Z_i)}\\
 &\qquad\le u_s\Bigl(d_1\normF{\L}\normF{\wh X_i^+-\wh X_i^-}
   +d_2\normF{\wh{\mathcal R}^+(\wh Z_i)-\wh{\mathcal R}^-(\wh Z_i)}\Bigr).
 \end{aligned}
\end{equation}
The matrices $\wh{\mathcal R}^\pm$ are reconstructed exactly from the factors supplied to the solver, and the constants $d_1,d_2$ in assumption~\ref{assump:ldlt-solver-smallness} are those in~\eqref{eq:assump-chol}.
It is reasonable to consider the rounding errors in solving the two correction equations jointly, as they are best solved together because they share the same coefficient matrix $A$ (see Section~\ref{sect:mp-ir-cost} below).

Since $P_i$ is a permutation and $J_i$ is a diagonal sign matrix, multiplication by either is exact. Hence, forming \(\widehat H_i\) in Fragment~\ref{alg.res-fac-chol} and \(\widehat K_i\) in Fragment~\ref{alg.sol-update-chol} requires one rounded matrix multiplication in each case.
Moreover,
\[
 \normt{P_i}=\normt{J_i}=1,\qquad
 \normF{P_i}=\sqrt{m_i},\qquad \normF{J_i}=\sqrt{p_i}.
\]
Under the same definitions of $b_2$ and $b_3$ as in~\eqref{eq:assump-c2-ldlt} and~\eqref{eq:assump-c3-ldlt}, we arrive at
\[
 \begin{aligned}
 \normF{T_i}^2
 &\le\sqrt{r_i}\normF{T_iT_i^T}
  =b_2\sqrt{r_i/m_i}\normF{H_i},\\
 \normF{\Gamma_i}^2
 &\le\sqrt{q_i}\normF{\Gamma_i\Gamma_i^T}
  =b_3\sqrt{q_i/p_i}\normF{K_i}.
 \end{aligned}
\]
For the reconstructed residual, we use~\eqref{eq:res-perturb-bound} with $\normt{P_i}=1$ in place of $\normF{N_i}$ and one rounded kernel product in forming $\widehat H_i$.
For the solution update, the analogous estimate follows from the backward QR relation~\eqref{eq:qr-Gi-hat} and the definition of $\wh K_i$ immediately before~\eqref{eq:deltaK-def-ldlt}, using $\normt{J_i}=1$ in place of $\normF{\Upsilon_i}$.
Substituting the preceding bounds for \(\|T_i\|_F^2\) and \(\|\Gamma_i\|_F^2\) into these perturbation estimates then gives
\begin{align*}
 \normF{\wh U_i\wh H_i\wh U_i^T-F_iP_iF_i^T}
 &\lesssim b_2\sqrt{r_i/m_i}
 \bigl(2(b_1+m_i^{3/2})\gammat_n^r+\gammat_{m_i}^r\bigr)\normF{H_i},\\
 \normF{\wh V_i\wh K_i\wh V_i^T-G_iJ_iG_i^T}
 &\lesssim b_3\sqrt{q_i/p_i}
 \bigl(2p_i^{3/2}\gammat_n^c+\gammat_{p_i}^c\bigr)\normF{K_i}.
\end{align*}

The Cholesky factors in Fragments~\ref{alg.res-fac-chol} and~\ref{alg.sol-update-chol} additionally contain square roots of the magnitudes of the retained eigenvalues.
For the residual factors, we adapt~\eqref{eq:residual-factor-weighted-bound-ldlt} by weighting each column of the factor error by the square root of the magnitude of the corresponding retained eigenvalue.
Applying the same argument to the solution factor in~\eqref{eq:Xhat-i+1-ldlt} yields the following first-order bounds:
\[
 \begin{aligned}
 \normF{\Delta R_i}&\lesssim
 2\sqrt{m_i}(u_r+2\gamma_{m_i}^r)\normF{\wh H_i},\\
 \normF{\Delta\wh X_{i+1}}&\lesssim
 2\sqrt{p_i}(u_c+2\gamma_{p_i}^c)\normF{\wh X_{i+1}}.
 \end{aligned}
\]
Adapting the column weighting in~\eqref{eq:residual-factor-weighted-bound-ldlt} to the combined retained positive and negative eigenpairs gives
\[
 \normF{\Delta R_i^s}\lesssim
 2u_s\normF{\wh{\mathcal R}^+(\wh Z_i)-\wh{\mathcal R}^-(\wh Z_i)},
\]
and thus $2u_s$ replaces $3u_s$ in~\eqref{eq:res-conversion-bound}.

Retaining first-order terms and neglecting the smaller rounding terms with indices $m_i$ and $p_i$, under assumption~\ref{assump:ldlt-qr}, we obtain the following variants of~\eqref{eq:phi-ldlt} and~\eqref{eq:sol-upt-delta-compute-normbnd},
\begin{equation}\label{eq:chol-error-coefficients}
 \begin{aligned}
 \phi_i&:=\psi+2\sqrt{m_i}\eta_r+2u_s
       +2b_2\sqrt{r_i/m_i}(b_1+m_i^{3/2})\gammat_n^r,\\
 \normF{\Delta\Xi_i}
 &\lesssim\Bigl[2\sqrt{p_i}\eta_s+\eta_n
       +2b_3\sqrt{q_i/p_i}\,p_i^{3/2}\gammat_n^c\Bigr]
       \normF{\wh X_{i+1}},
 \end{aligned}
\end{equation}
where $\psi=(d_1\kappa_F(\L)+d_2)u_s$ uses the constants in~\eqref{eq:assump-chol}.
The recurrence~\eqref{eq:common-residual-recurrence} and Theorem~\ref{thm:res-converge-ldlt} therefore apply with the coefficients in~\eqref{eq:chol-error-coefficients} and $\wh X_{i+1}=\wh Z_{i+1}\wh Z_{i+1}^T$, provided $\phi:=\sup_{i\ge1}\phi_i<1$ and the corresponding supremum in the bound for the limiting residual is finite.

\section{The sign function Newton iteration}\label{sect.sign-func-Newton-iter}
The matrix sign function method for solving the Lyapunov equation was introduced by Roberts in 1971~\cite{robe80}, and it has since been one of the most widely used methods, owing to its easy yet robust implementation, excellent parallelism, and richness in level-3 BLAS operations~\cite{hqsw00}.
To test numerically the convergence and attainable accuracy predicted by the rounding error analysis of the mixed-precision IR framework, we use the \emph{factorized} sign function Newton iteration as a dense test solver.
Its factored formulations directly provide the Cholesky-type and $\ldlt$-type factors required by the IR algorithms.
The IR framework and its error analysis also allow different solvers, such as low-rank ADI~\cite{blp08,lms15} or Krylov methods~\cite{jbri06,simo07}.

The (scaled) sign function Newton iteration for the solution $X$ of Lyapunov equations~\eqref{eq:lyap} is
\begin{subequations}\label{eq:newton-lyapu}
	\begin{align}
		A_{k} &= \frac{1}{2} \bigl( \mukm A_{k-1} + \mukm^{-1}A_{k-1}^{-1}\bigr),
		&  A_0 &= A,  \label{eq:newton-lyapu1} \\
		W_{k} &= \frac{1}{2} \bigl( \mukm W_{k-1} +\mukm^{-1} A_{k-1}^{-1}W_{k-1}A_{k-1}^{-T}\bigr),
		&  W_0 &= W,  \label{eq:newton-lyapu2}
	\end{align}
\end{subequations}
where the scaling parameter $\mukm>0$ can be used to
accelerate the convergence of the method in its initial steps. 
The choice $\muk \equiv 1$ yields the unscaled Newton iteration, for which $A_k$ and $W_k$ converge quadratically to $-I_n$ and $2X$, respectively~\cite{bequ99}, \cite{robe80}.  
Common scaling techniques include the determinantal scaling $\mu_k = |\det(A_k)|^{-1/n}$, the spectral scaling $\mu_k = \rho(A_k^{-1})^{1/2} / \rho(A_k)^{1/2}$, and the $2$-norm scaling $\mu_k = \normt{A_k^{-1}}^{1/2} / \normt{A_k}^{1/2}$. The spectral norm is often approximated by the Frobenius norm when it is expensive to calculate~\cite{high86}, so there is also the Frobenius-norm scaling $\mu_k = \normF{A_k^{-1}}^{1/2} / \normF{A_k}^{1/2}$, 
since the computation of the Frobenius norm parallelizes quite well---each matrix entry contributes independently to the final result. 
In general, there is no single scaling strategy that is superior to the rest~\cite[sect.~5.5]{high08}, \cite{sibe08}.

\subsection{Iterating on the solution factors}
In the case of $W=LL^T$,
Larin and Aliev~\cite{laal93} propose to rewrite the iteration for $W_k$ in~\eqref{eq:newton-lyapu2} as
\begin{equation*}
	W_{k} = \frac{\mukm}{2} 
	\begin{bmatrix}
		Z_{k-1} & \mukm^{-1}A_{k-1}^{-1}Z_{k-1}
	\end{bmatrix}
	\begin{bmatrix}
		Z_{k-1}^T \\ \mukm^{-1}Z_{k-1}^TA_{k-1}^{-T}
	\end{bmatrix},
	\quad  W_0 = Z_0Z_0^T \equiv LL^T,
\end{equation*}
and thus obtain the iterations for the Cholesky-type factored solution
\begin{subequations}\label{eq:newton-lyapu-factored-chol}
	\begin{align}
		A_{k} &= \frac{1}{2} \bigl( \mukm A_{k-1} + \mukm^{-1}A_{k-1}^{-1}\bigr),
		\hspace{-80pt } &  A_0 &= A,   \label{eq:newton-lyapu-factored-chol1}\\
		Z_{k} &= \sqrt{\frac{\mukm}{2}} 
		\begin{bmatrix}
			Z_{k-1} & \mukm^{-1}A_{k-1}^{-1}Z_{k-1}
		\end{bmatrix},
		\hspace{-80pt } &  Z_0 &= L,  \label{eq:newton-lyapu-factored-chol2}
	\end{align}
\end{subequations}
where $Z_kZ_k^T/2$ converges to $X$ without rank truncation.

Writing $W=LSL^T$ for a symmetric matrix $S$ (positive semidefinite for the original equation, possibly indefinite for a correction equation), we can obtain the iterations for the solution factors
\begin{subequations}\label{eq:newton-lyapu-factored-ldlt}
	\begin{align}
		A_{k} &= \frac{1}{2} \bigl( \mukm A_{k-1} + \mukm^{-1}A_{k-1}^{-1}\bigr),
		\hspace{-80pt } &  A_0 &= A,   \label{eq:newton-lyapu-factored-ldlt1}\\
		Y_{k} &=  
		\begin{bmatrix}
			\frac{\mukm}{2} Y_{k-1} &  \\
			&  \frac{1}{2\mukm} Y_{k-1}
		\end{bmatrix},
		\hspace{-80pt } &  Y_0 &= S, 
		\label{eq:newton-lyapu-factored-ldlt2}\\
		Z_{k} &= 
		\begin{bmatrix}
			Z_{k-1} & A_{k-1}^{-1}Z_{k-1}
		\end{bmatrix},
		\hspace{-80pt } &  Z_0 &= L, 
		\label{eq:newton-lyapu-factored-ldlt3}
	\end{align}
\end{subequations}
where $Z_kY_kZ_k^T/2$ converges to $X$ without rank truncation.
Note that this $\ldlt$-type formulation avoids scaling on the tall-and-skinny solution factor $Z$, which could be expensive when the columns of $Z$ accumulate as the number of iterations grows. 
The iteration on the inner factor $Y_k$ can be reduced to one involving only the scaling parameters $\muk$, and $Y_k$ can be formed at the end as
\begin{equation*}Y_{k} =  
	\begin{bmatrix}
		\frac{\mukm}{2} &  \\
		&  \frac{1}{2\mukm} 
	\end{bmatrix}\otimes \cdots \otimes
	\begin{bmatrix}
		\frac{\mu_1}{2} &  \\
		&  \frac{1}{2\mu_1} 
	\end{bmatrix}\otimes 
	\begin{bmatrix}
		\frac{\mu_0}{2} &  \\
		&  \frac{1}{2\mu_0} 
	\end{bmatrix}\otimes S.
\end{equation*}
The $\ldlt$-type formulation has essentially the same storage requirement as the Cholesky-type formulation.

\subsection{The solvers}
The factorized iterations~\eqref{eq:newton-lyapu-factored-chol}--\eqref{eq:newton-lyapu-factored-ldlt} can substantially reduce the computational cost, as it avoids the two full $n$-by-$n$ matrix multiplications required  in~\eqref{eq:newton-lyapu}. Nevertheless, a potential concern is the growth of the size of the iterate $Z_k$, whose column dimension $c_k$ doubles at each iteration, which implies that the required storage space grows exponentially. 
Therefore, low-rank truncations are often conditionally performed within the iteration to limit the increase of the smaller dimension of the low-rank factors~\cite{bekq11},~\cite{bequ99},~\cite{laal93},~\cite[sect.~7.1]{schu22}.

\begin{algorithm2e}[t]
	\caption{The sign function Newton iteration for Cholesky-type solution factor of the low-rank Lyapunov equation~\eqref{eq:lyap}.}
	\label{alg.newton-lyapu-chol}
	\SetAlgoVlined
	\Parameter{Unit roundoff $u>0$, convergence tolerance $\tau_N=10\sqrt{nu}$, maximal iterations $k_{\max}\in\mathbb{N}^+$, rank truncation threshold $0<\vartheta\le 1$}
	{
		\KwIn{$A\in\R^{\nbyn}$, $L\in\R^{n\times m}$}
		\KwOut{Approximated solution factor $Z$ such that $X\approx ZZ^T$}
		$Z_0= L$, $A_0= A$\;
		\For{$k \gets 1$ \KwTo $k_{\max}$}{
			$A_k = \frac{1}{2} ( \mukm A_{k-1} + \mukm^{-1}A_{k-1}^{-1} )$ with $\mukm$ a scaling parameter. \;
			$Z_{k} =  \frac{1}{\sqrt{2}} \begin{bmatrix}
				\mukm^{1/2} Z_{k-1} &  
				\mukm^{-1/2} A_{k-1}^{-1}Z_{k-1}
			\end{bmatrix}$ \;
			\label{algline.newton-lyapu-matmul}
			\lIf{$\size(Z_k, 2) > \vartheta n$}{
				$Z_k$ = \textsc{RankTruncChol}$(Z_k)$  \label{algline.newton-lyapu-ranktrunc}}
			\If{$\normi{A_k + I_n} \le \tau_N$}{Terminate after one further iteration. \; }
		}
		$Z = Z_{k} / \sqrt{2}$ \;
		\SetAlgoNoLine
	}
	\lineseparator
	\nonl
	\subFn{$Z$ = \textsc{RankTruncChol}$(Z)$
	}{
		Compute a rank-revealing QR:
		$Z^T\mPi = Q\begin{bmatrix}
			T & C \\ 0 & S
		\end{bmatrix}$ with $\normt{S}\le\sqrt{u}\normt{Z}$. \;
\Return{$\mPi\begin{bmatrix}
				T  & C
			\end{bmatrix}^T$} 
	}
\end{algorithm2e}

\begin{algorithm2e}[t]
	\caption{The sign function Newton iteration for $\ldlt$-type solution factor of the low-rank Lyapunov equation~\eqref{eq:lyap}.}
	\label{alg.newton-lyapu-ldlt}
	\SetAlgoVlined
	\Parameter{Unit roundoff $u>0$, convergence tolerance $\tau_N=10\sqrt{nu}$, maximal iterations $k_{\max}\in\mathbb{N}^+$, rank truncation threshold $0<\vartheta\le 1$}
	{
		\KwIn{$A\in\R^{\nbyn}$, $L\in\R^{n\times m}$, and $S\in\R^{m\times m}$}
		\KwOut{Approximated solution factors $Z$ and $Y$ such that $X\approx ZYZ^T$}
		$Z_0= L$, $Y_0= S$, $A_0= A$\;
		\For{$k \gets 1$ \KwTo $k_{\max}$ }{
			$A_k = \frac{1}{2} ( \mukm A_{k-1} + \mukm^{-1}A_{k-1}^{-1} )$ with $\mukm$ a form of scaling. \;
			$Z_k = \begin{bmatrix}
				Z_{k-1} &   A_{k-1}^{-1} Z_{k-1}
			\end{bmatrix}$ \;
			\label{algline.newton-lyapu-matmul-ldlt}
			$Y_k= \frac{1}{2} \blkdiag(\mukm Y_{k-1},\ \mukm^{-1}Y_{k-1}) $ \;
			\lIf{$\size(Z_k, 2) > \vartheta n$}{
				$[Z_k, Y_k]$ = \textsc{RankTruncLdlt}$(Z_k, Y_k)$ 
				\label{algline.newton-lyapu-ranktrunc-ldlt}}
			\If{$\normi{A_k + I_n} \le \tau_N$}{Terminate after one further iteration. \; }
		}
		$Z=Z_k$, $Y = Y_k / 2$ \;
		\SetAlgoNoLine
	}
	\lineseparator
	\nonl
	\subFn{$[Z, Y]$ = \textsc{RankTruncLdlt}$(Z, Y)$
	}{
		Compute a thin QR decomposition $Z = QR$. \;
		Compute the spectral decomposition $RYR^T = \wt{V}\wt{\Lambda} \wt{V}^T$. \;
		$\Lambda = \diag(\wt{\lambda}_i)$, 
		$i\in I_\lambda := \{i\ \vert\  |\wt{\lambda}_i| > u \normi{ \wt{\Lambda}} \}$ and 
		$V= \wt{V}(\colon, I_\lambda)$ \;
\Return{$QV$, $\Lambda$} 
	}
\end{algorithm2e}

The overall algorithms for both types of solvers are presented as Algorithm~\ref{alg.newton-lyapu-chol} and Algorithm~\ref{alg.newton-lyapu-ldlt}. 
The stopping criteria for both iterations depend only on the coefficient matrix $A$, which is set as 
\begin{equation}\label{eq:newton-stop-tau}
	\normi{A_k + I_n} \le \tau_N,\quad \tau_N = 10\sqrt{nu},
\end{equation}
with one additional iteration being performed after the tolerance has been reached.
This is to avoid potential stagnation, which occurs when the stopping 
criterion is too stringent, yet still letting the algorithm try to reach the attainable accuracy; 
see~\cite{bequ99} for more details on the setting of the convergence test.

\subsection{Computational cost analysis of the IR algorithms}
\label{sect:mp-ir-cost}
We discuss the cost of Algorithms~\ref{alg.ir3-lyap-chol} and~\ref{alg.ir3-lyap-ldlt} when the inner solver is the sign function Newton iteration, Algorithm~\ref{alg.newton-lyapu-chol} or~\ref{alg.newton-lyapu-ldlt}.

Recall that the coefficient matrix $L$ has dimensions $m\ll n$ and the sought solution factor $Z$ is expected to have small numerical rank with respect to $n$, so one can safely choose a $\vartheta\ll 1$ for deciding when to perform a rank truncation.
Ideally, the rank truncation threshold should be chosen such that $\vartheta n$ is greater than $m$ and the numerical rank of $Z$, the latter of which is, however, not known a priori.

As a baseline, we can assume that the total number of rank truncations $t$ is in general much smaller than the total number of iterations $k$ for convergence, i.e., $t\ll k$.
Since the possible rank truncation imposes a constraint on the size of the iterates, in any iteration,
the smaller matrix in the matrix product in Line~\ref{algline.newton-lyapu-matmul} of Algorithm~\ref{alg.newton-lyapu-chol} or Line~\ref{algline.newton-lyapu-matmul-ldlt} of Algorithm~\ref{alg.newton-lyapu-ldlt} is at most of size $n\times \vartheta n$ and 
the rank truncation in Line~\ref{algline.newton-lyapu-ranktrunc} of Algorithm~\ref{alg.newton-lyapu-chol} or Line~\ref{algline.newton-lyapu-ranktrunc-ldlt} of Algorithm~\ref{alg.newton-lyapu-ldlt} is performed on a matrix no larger than $n\times 2\vartheta n$.
It follows that, in each iteration, the matrix product requires $2\vartheta n^3$ flops and the rank-revealing QR costs $O(\vartheta^2n^3)$ flops; see~\cite[sect.~5.4.2]{gova13} for the flops count for Householder QR with column pivoting, for example.
Overall, Algorithm~\ref{alg.newton-lyapu-chol} and Algorithm~\ref{alg.newton-lyapu-ldlt} require 
$2k n^3 + O(\vartheta k n^3)$ flops for $k$ Newton iterations performed.

The Cholesky-type formulation in each refinement step solves two correction equations~\eqref{eq:lyap-cor-chol}, which share the common coefficient matrix $A$.
Therefore, one can easily adapt Lines~\ref{algline.newton-lyapu-matmul}--\ref{algline.newton-lyapu-ranktrunc} of Algorithm~\ref{alg.newton-lyapu-chol} to solve the two equations simultaneously, such that the adapted algorithm produces the two solution factors $Z_i^+$ and $Z_i^-$ together. 
This saves a full $n\times n$ matrix inversion per iteration but preserves the ability to perform the rank truncations of the factors independently.
Therefore, the asymptotic cost of invoking Algorithm~\ref{alg.newton-lyapu-chol} for solving the correction equations~\eqref{eq:lyap-cor-chol} remains $2 k n^3 + O(\vartheta k n^3)$ flops in total.

For Fragment~\ref{alg.res-fac-chol}--Fragment~\ref{alg.sol-update-chol} and Fragment~\ref{alg.res-fac-ldlt}--Fragment~\ref{alg.sol-update-ldlt},
the thin QR decompositions are performed on a tall-and-skinny matrix (with the smaller dimension much lower than $n$).
As a result, the cost of a single call to one of these subroutines is only $O(n^2)$ flops, which is negligible compared with an invocation of the solver, under the practical condition that a flop at precision $u_c$ or $u_r$ is not much more expensive than $n$ flops at precision $u_s$.

To sum up, when the sign function Newton iterations, Algorithm~\ref{alg.newton-lyapu-chol} and Algorithm~\ref{alg.newton-lyapu-ldlt}, are used as the solver,
the asymptotic cost of the IR scheme with $i$ refinement steps performed is $2(k_0+k_1+\cdots+k_i)n^3$ flops at precision $u_s$, where $k_0$ denotes the number of Newton steps for the initial solve and $k_\ell$, $\ell=1\colon i$, denotes the number of \textit{inner} Newton steps carried out at the $\ell$th \textit{outer} refinement step.
The dominant cost of the algorithm therefore depends on both the precision of the solver and the number of total Newton iterations performed throughout. 
When a lower precision for $u_s$ is used, the convergence of the Newton iteration is expected to be slower, leading to a higher number of iterations. This is reflected by the varying condition bounds for different combinations of $u_s$, $u_c$, and $u_r$ in Table~\ref{table:ir3}. 

\subsubsection{Computational cost model}
\label{sect:mp-ir-cost-model}
To model the effect of the solver precision on the computational cost, we fix the working precision, the precisions for the residual factorization and solution update and compare a lower solver precision with a higher one.
For a fixed problem size, let $t_{\mathrm{low}}>0$ and $t_{\mathrm{high}}>0$ denote effective times per flop for the dominant matrix-inversion operation at the two solver precisions.
Following the flop-based cost model in~\cite[sect.~5]{dfhl25}, we parameterize their relative cost by the ratio
\begin{equation}\label{eq:rho-cost}
 \rho:=\frac{t_{\mathrm{low}}}{t_{\mathrm{high}}}>0.
\end{equation}
Computing at higher precision is usually more expensive per flop, so we expect $\rho\le1$ for implementations that efficiently exploit lower precision, and potentially $\rho\ll1$ on emerging hardware with dedicated support for low-precision arithmetic; cf.~Section~\ref{sect.intro}.

For each run, the total number of inner Newton steps over the initial solve and all $i$ outer refinement steps is
\begin{equation}\label{eq:cost-iteration-counts}
 k^{\sigma}:=\sum_{\ell=0}^{i}k_\ell.
\end{equation}
Using the leading inversion cost of $2n^3$ flops per Newton step, we model the dominant contribution to the computational time as
\[
 T_{\mathrm{low}}^{\mathrm{inv}}
 =2n^3t_{\mathrm{high}}\rho k^{\sigma}_{\mathrm{low}},
 \qquad
 T_{\mathrm{high}}^{\mathrm{inv}}
 =2n^3t_{\mathrm{high}}k^{\sigma}_{\mathrm{high}}.
\]
We define the conditional speedup as the modeled time at the higher solver precision divided by that at the lower solver precision:
\begin{equation}\label{eq:speedup-no-reuse}
 \mathcal{S}_{\sigma}(\rho)
 :=\frac{T_{\mathrm{high}}^{\mathrm{inv}}}{T_{\mathrm{low}}^{\mathrm{inv}}}
 =\frac{k^{\sigma}_{\mathrm{high}}}{\rho k^{\sigma}_{\mathrm{low}}},
 \qquad
 \mathcal{S}_{\sigma}(\rho)>1
 \ \Longleftrightarrow\ 
 \rho<\rho_{\sigma}^{*}:=\frac{k^{\sigma}_{\mathrm{high}}}{k^{\sigma}_{\mathrm{low}}}.
\end{equation}
This speedup is conditional on the assumed value of $\rho$ and on both runs attaining the prescribed residual tolerance, although their numbers of inner Newton steps and outer refinements may differ; it accounts only for the dominant inversion work.
At the break-even threshold $\rho=\rho_{\sigma}^{*}$, the two modeled dominant computational times are equal.
The lower solver precision reduces the dominant computational time when $\rho<\rho_{\sigma}^{*}$ and increases it when $\rho>\rho_{\sigma}^{*}$.
For example, if the lower-precision run requires twice as many inner Newton steps in total, then $\rho_{\sigma}^{*}=1/2$, so its inversion cost per flop must be less than half that at the higher precision to yield a speedup.

\subsubsection{Alternative cost model with reuse of the coefficient iterates}
\label{sect:mp-ir-cost-alt}
The correction equations~\eqref{eq:lyap-cor-chol} and~\eqref{eq:lyap-cor-ldlt} share the same coefficient matrix $A$.
For a fixed solver precision and an identical, deterministic scaling strategy, its Newton coefficient iterates, scaling factors, and computed inverses are therefore independent of the right-hand side.
They can be stored and reused across the outer refinement steps, extending the stored sequence only when a call requires additional Newton steps.
Let 
\begin{equation}\label{eq:cost-iteration-counts-max}
k^{\max}:=\max_{0\leq\ell\leq i}k_\ell
\end{equation}
denote the largest number of Newton steps required by an individual solver call.
With this reuse, only $k^{\max}$ matrix inversions are computed over the entire IR algorithm, including the initial solve and all $i$ outer refinement steps.
Their leading cost is $2n^3k^{\max}$ flops at precision $u_s$, giving the analogue of~\eqref{eq:speedup-no-reuse}:
\begin{equation}\label{eq:speedup-reuse}
 \mathcal{S}_{\max}(\rho)
 :=\frac{k^{\max}_{\mathrm{high}}}{\rho k^{\max}_{\mathrm{low}}},
 \qquad
 \rho_{\max}^{*}:=\frac{k^{\max}_{\mathrm{high}}}{k^{\max}_{\mathrm{low}}}.
\end{equation}
Note that each inner Newton step in Algorithms~\ref{alg.newton-lyapu-chol} and~\ref{alg.newton-lyapu-ldlt} must still form $A_{k-1}^{-1}Z_{k-1}$, which depends on the right-hand side of the current solver call.
Any required rank truncations must also be performed separately on the factors of each solve, including the QR factorizations and, in the $\ldlt$-type formulation, the small spectral decompositions.
The total cost of these products and truncations over the initial solve and all outer refinement steps therefore continues to depend on the total inner iteration count $k^{\sigma}$.
Under the same rank bounds as in the preceding cost analysis, this cost is $O(\vartheta n^3k^{\sigma})$ flops at precision $u_s$.
This contribution remains non-dominant relative to the inversion cost $2n^3k^{\max}$ provided that $\vartheta k^{\sigma}\ll k^{\max}$.

\begin{table}[t]
	\centering
	\caption{Approximate storage of the dense inverse sequence in native real formats. Other arrays and temporary workspace are excluded. Here $1\,\mathrm{MiB}=2^{20}$ bytes and $1\,\mathrm{GiB}=2^{30}$ bytes.}
	\label{tab:inverse-sequence-storage}
	\begin{tabular}{rrrrr}
		\toprule
		$n$ & $k^{\max}$ & bf16/fp16 & fp32 & fp64 \\
		\midrule
		$500$ & $10$ & $4.8$ MiB & $9.5$ MiB & $19.1$ MiB \\
		$1{,}000$ & $20$ & $38.1$ MiB & $76.3$ MiB & $152.6$ MiB \\
		$5{,}000$ & $20$ & $953.7$ MiB & $1.86$ GiB & $3.73$ GiB \\
		$10{,}000$ & $50$ & $9.31$ GiB & $18.63$ GiB & $37.25$ GiB \\
		\bottomrule
	\end{tabular}
\end{table}

Retaining the $k^{\max}$ dense inverse matrices requires $M_{\mathrm{seq}}=b n^2k^{\max}$ bytes, where $b$ is the number of bytes per stored entry.
For native bf16/fp16, fp32, and fp64 storage, $b=2$, $4$, and $8$, respectively; Table~\ref{tab:inverse-sequence-storage} reports the memory required to store the inverse sequence for selected matrix dimensions \(n\), sequence lengths \(k^{\max}\), and floating-point formats.
Keeping this sequence resident in cache reduces traffic to RAM, provided that sufficient cache capacity remains alongside the other active arrays.
Otherwise, repeated reads from RAM incur additional data-movement costs.
For a fixed sequence length, lower-precision storage reduces both the memory requirement and the number of bytes transferred per read, which can mitigate this overhead.
When comparing runs at different precisions, the total storage also depends on their respective values of $k^{\max}$.
The speedup in~\eqref{eq:speedup-reuse} therefore accounts only for the dominant inversion work and best approximates the full speedup when the factor operations and data-movement costs are small.

\def\Y{\widehat{Y}}
\def\Z{\widehat{Z}}
\def\res{{\text{res}}}
\section{Numerical experiments}\label{sect.experiments}
In this section, we test the attainable accuracy and convergence behaviour of the mixed-precision IR frameworks using the sign function Newton iteration, Algorithm~\ref{alg.newton-lyapu-chol} or~\ref{alg.newton-lyapu-ldlt}, as the inner solver.
We gauge the quality of computed solution factors by the relative residual of the equation~\eqref{eq:lyap} in the Frobenius norm, given by
\begin{equation}\label{eq:res-fac-sol}
\res(\Z) = \frac{\normF{A\Z\Z^T+\Z\Z^T A^T+W}}{\normF{W} + 2\normF{\Z\Z^T}\normF{A}},
\quad \mbox{or} \quad
\res(\Z, \Y) = \frac{\normF{A\Z\Y\Z^T+\Z\Y\Z^T A^T+W}}{\normF{W} + 2\normF{\Z\Y\Z^T}\normF{A}}, 
\end{equation} 
depending on the type of the solver used in the algorithm.
All reported residuals are evaluated at termination from the returned solution factors.
For the tests with fixed solver-precision, matrix products and norms in~\eqref{eq:res-fac-sol} are computed in the working precision; in the adaptive solver-precision tests, they are computed in fp64.
The IR stopping test uses residual norms estimated from all eigenvalues of the small symmetric residual cores in Fragments~\ref{alg.res-fac-chol} and~\ref{alg.res-fac-ldlt}, before rank truncation.

All experiments were performed on a machine equipped with a 2.0 GHz quad-core Intel Core i5-1038NG7 processor and 16 GiB of RAM, using four computational threads in MATLAB 25.2 (R2025b Update 5, 64-bit macOS). The code that produces the results in this section is available on GitHub.\footnote{\url{https://github.com/xiaobo-liu/mpir-lrlyap}}
We use bf16 as the half precision format, with rounding simulated using the \texttt{chop}\footnote{\url{https://github.com/higham/chop}} function~\cite{hipr19}.
The experiments therefore aim to validate the error analysis and supply iteration counts for the conditional cost models~\eqref{eq:speedup-no-reuse} and~\eqref{eq:speedup-reuse}.

We tried the different scaling schemes mentioned in Section~\ref{sect.sign-func-Newton-iter} as well as the combined scaling (with the Frobenius norm and the determinant) recommended in~\cite{sibe08} 
for the Newton iteration, but we found no technique bringing benefits dominating the others for our test sets. 
In particular, scalings that require the computation of the determinant of a large matrix are more prone to suffer from underflow and overflow issues in low precision.
In our implementation we use the Frobenius-norm scaling, where we also monitor the relative change of the iterates,
\begin{equation*}
	\delta_k := \normF{A_k - A_{k-1}} / \normF{A_k},
\end{equation*}
for the use of scaling and premature termination of the iterations. 
We adopt the strategy of Higham~\cite[sect.~5.8]{high08} and stop scaling of the iterations once $\delta_k< 10^{-2}$, which is to avoid the interference of nonoptimal scaling parameters on the convergence when the region of convergence is reached. 
Higham also found that the condition
\begin{equation}\label{eq:newton-inner-stagnation}
\delta_k>\delta_{k-1}/2,
\end{equation}
which means that $\delta_k$ has not decreased by at least a factor $2$, is a good indicator of the dominance of roundoff errors.
In the implementation, either this condition or the prescribed test~\eqref{eq:newton-stop-tau} triggers termination after one further Newton step, subject to the iteration limit; all reported iteration counts follow this rule.

In our implementation of Algorithm~\ref{alg.ir3-lyap-chol} and Algorithm~\ref{alg.ir3-lyap-ldlt}, we monitor the ratio of two successive relative residuals~\eqref{eq:res-fac-sol}, defined by  
\begin{equation}\label{eq:it-early-terminate}
	\theta_i := \res_i / \res_{i-1}.
\end{equation}
A ratio $\theta_i$ close to $1$ means little improvement has been made in the previous refinement step. We therefore terminate the refinement process if $\theta_i>0.9$ (the residual is reduced by less than $10\%$) for two consecutive iterations. We found this scheme can effectively signify stagnation of the refinement process, especially for ill-conditioned problems.

\subsection{Specification of the algorithmic parameters}\label{sect:algorithmic-parameters}

The global convergence tolerance of the IR framework is set to $\tau_I=nu$, with a maximum of $i_{\max} = 50$ refinement steps.
The maximal iteration for Algorithm~\ref{alg.newton-lyapu-chol} and Algorithm~\ref{alg.newton-lyapu-ldlt} is set to $k_{\max}=50$.
The rank truncation tolerance $\vartheta$ controls when the factors in the Newton solver are compressed, but anoptimal, or even suitable, value of $\vartheta$ depends on the numerical rank of the sought solution factor, as discussed in Section~\ref{sect:mp-ir-cost}.
Since the coefficient matrix $L$ in~\eqref{eq:lyap}, in our experiments, has dimension $m< 0.1 n$, we set $\vartheta=0.1$ correspondingly.

According to the analysis in Section~\ref{sect:err-anal-ir-ldlt} and its Cholesky specialization,
the spectral splitting tolerance $\eta_s>0$ in Fragment~\ref{alg.sol-update-chol} and Fragment~\ref{alg.sol-update-ldlt} is set to $10u$, 
and we choose the residual spectral splitting threshold $\eta_r=10^{-4}$ in Fragment~\ref{alg.res-fac-chol} and Fragment~\ref{alg.res-fac-ldlt}.

\subsection{Tests on synthetic matrices}
The experiments of this section are on low-rank Lyapunov equations constructed using pseudorandom matrices with specified order of condition number. 
The coefficient matrices in~\eqref{eq:lyap} were generated with the MATLAB code
\begin{lstlisting}
	L = rand(n, m);
	W = L * L.'; 
	V = gallery('orthog', n);
	A = - V .* (logspace(0, q, n)) * V.';
\end{lstlisting}
setting $m=3$. 
Note that this code generates a symmetric coefficient matrix $A$, but the symmetry is not explicitly exploited in the algorithms; the construction of $A$ is for controlling the condition of the Lyapunov equation, such that $\kappa_F(\L)\approx 10^q$, $q\ge 0$.
For the $\ldlt$-type solver, the inner factor matrix $S$ of $W=LSL^T$ is initialized to be the identity matrix of order $m$.

We examine the quality of the computed solutions by~\eqref{eq:res-fac-sol} under the precision settings listed in Table~\ref{table:ir3}.
The sizes of the coefficient matrix $A$ are set to $n=100$ and $n=1000$.
In total, with varying condition numbers and sizes, 16 different low-rank Lyapunov equations~\eqref{eq:lyap} tested.

\newcommand{\iter}{\text{iter}}
\newcommand{\rk}{\text{rank}}
\begin{table}[t]
	\renewcommand{\tabcolsep}{4pt}
	\centering
	\caption{Results on low-rank Lyapunov equations of varying sizes and condition numbers. For each solver precision $u_s$, 
		iter = $\varsigma(\alpha)$ means the total number of Newton iterations over all refinement steps is $\varsigma$ and the maximal number of iterations in a call is $\alpha$. A dash ``--'' in the rank indicates failure to converge to the prescribed residual tolerance.}
	\resizebox{1.00\linewidth}{!}{
		\begin{tabularx}{1.3\textwidth}{@{\extracolsep{\fill}}cc|lrr|lrr|lrr|lrr|lrr}
\toprule
\multicolumn{2}{c|}{Cholesky-type} & \multicolumn{6}{c|}{$u=$ fp32} & \multicolumn{9}{c}{$u=$ fp64} \\
\midrule
\multicolumn{2}{c|}{} & \multicolumn{3}{c|}{$u_s=$ bf16} & \multicolumn{3}{c|}{$u_s=$ fp32} & \multicolumn{3}{c|}{$u_s=$ bf16} & \multicolumn{3}{c|}{$u_s=$ fp32} & \multicolumn{3}{c}{$u_s=$ fp64} \\
$n$ & $\cond$ & $\res$ & $\iter$ & $\rk$ & $\res$ & $\iter$ & $\rk$ & $\res$ & $\iter$ & $\rk$ & $\res$ & $\iter$ & $\rk$ & $\res$ & $\iter$ & $\rk$ \\
\midrule
$100$ & 1.3e0 & 1.4e-6 &  4(2) &   6 & 6.7e-7 &  2(2) &   6 & 6.9e-16 & 16(2) &  15 & 6.8e-15 &  6(2) &  17 & 7.0e-17 &  4(4) &  15   \\
         & 3.2e0 & 3.2e-7 &  8(2) &   9 & 3.6e-7 &  3(3) &   9 & 3.5e-15 & 18(2) &  24 & 3.3e-15 &  9(3) &  24 & 7.1e-17 &  5(5) &  24   \\
         & 1.0e1 & 6.7e-7 &  8(2) &  12 & 4.4e-8 &  4(4) &  12 & 7.1e-16 & 22(2) &  33 & 2.2e-16 & 12(4) &  33 & 1.8e-15 &  5(5) &  33   \\
         & 3.2e1 & 2.0e-6 & 12(2) &  13 & 3.1e-7 &  4(4) &  13 & 7.9e-15 & 36(2) &  42 & 3.5e-15 & 12(4) &  45 & 9.9e-17 &  6(6) &  42   \\
         & 1.0e2 & 5.7e-6 & 20(2) &  20 & 5.2e-8 &  5(5) &  15 & 7.9e-15 & 78(2) &  49 & 8.0e-16 & 15(5) &  49 & 1.0e-16 &  6(6) &  49   \\
         & 3.2e2 & 3.0e-6 & 27(3) &  18 & 7.0e-8 &  5(5) &  16 & 6.5e-15 & 75(3) &  58 & 6.8e-16 & 15(5) &  57 & 7.8e-17 &  7(7) &  56   \\
         & 1.0e3 & 4.4e-4 & 24(3) & -- & 3.7e-8 &  5(5) &  16 & 2.9e-4 & 21(3) & -- & 3.8e-16 & 15(5) &  62 & 8.3e-17 &  7(7) &  62   \\
         & 3.2e3 & 2.5e-3 & 15(5) & -- & 6.2e-8 &  5(5) &  16 & 3.2e-3 & 15(5) & -- & 1.7e-15 & 15(5) &  71 & 7.4e-17 &  7(7) &  68   \\
\midrule
$1000$ & 1.3e0 & 2.9e-7 &  4(2) &   6 & 2.4e-7 &  2(2) &   6 & 1.0e-13 & 16(2) &  44 & 2.7e-16 &  6(2) &  12 & 1.9e-17 &  4(4) &  12   \\
         & 3.2e0 & 1.1e-6 &  4(2) &   6 & 1.3e-7 &  3(3) &   6 & 6.6e-14 & 20(2) &  75 & 1.6e-15 &  9(3) &  22 & 1.8e-17 &  5(5) &  22   \\
         & 1.0e1 & 1.1e-5 &  4(2) &   9 & 7.3e-6 &  3(3) &   9 & 8.5e-14 & 26(2) &  44 & 9.6e-15 & 12(3) &  31 & 8.6e-16 &  5(5) &  31   \\
         & 3.2e1 & 2.6e-5 &  6(2) &  10 & 1.2e-7 &  4(4) &  10 & 2.6e-5 & 10(2) & -- & 2.4e-15 & 12(4) &  39 & 2.8e-16 &  6(6) &  39   \\
         & 1.0e2 & 7.4e-5 & 10(2) & -- & 5.7e-7 &  4(4) &  11 & 7.4e-5 & 10(2) & -- & 1.6e-16 & 16(4) &  48 & 2.1e-16 &  6(6) &  48   \\
         & 3.2e2 & 2.2e-4 & 10(2) & -- & 7.4e-9 &  5(5) &  12 & 2.2e-4 & 10(2) & -- & 2.1e-16 & 15(5) &  55 & 3.4e-16 &  7(7) &  55   \\
         & 1.0e3 & 1.0e-4 &  9(3) & -- & 8.2e-9 &  5(5) &  12 & 1.1e-4 &  9(3) & -- & 1.6e-15 & 15(5) &  70 & 2.7e-16 &  7(7) &  61   \\
         & 3.2e3 & 5.7e-5 & 15(3) &  17 & 1.1e-8 &  5(5) &  10 & 9.6e-6 & 48(3) & -- & 2.4e-14 & 15(5) &  94 & 2.9e-16 &  7(7) &  69   \\
\midrule
\multicolumn{2}{c|}{$\ldlt$-type} & \multicolumn{6}{c|}{$u=$ fp32} & \multicolumn{9}{c}{$u=$ fp64} \\
\midrule
\multicolumn{2}{c|}{} & \multicolumn{3}{c|}{$u_s=$ bf16} & \multicolumn{3}{c|}{$u_s=$ fp32} & \multicolumn{3}{c|}{$u_s=$ bf16} & \multicolumn{3}{c|}{$u_s=$ fp32} & \multicolumn{3}{c}{$u_s=$ fp64} \\
$n$ & $\cond$ & $\res$ & $\iter$ & $\rk$ & $\res$ & $\iter$ & $\rk$ & $\res$ & $\iter$ & $\rk$ & $\res$ & $\iter$ & $\rk$ & $\res$ & $\iter$ & $\rk$ \\
\midrule
$100$ & 1.3e0 & 1.6e-6 &  4(2) &   6 & 6.7e-7 &  2(2) &   6 & 1.2e-15 & 12(2) &  15 & 4.1e-15 &  6(2) &  17 & 1.5e-16 &  4(4) &  15   \\
         & 3.2e0 & 3.9e-6 &  4(2) &  10 & 3.5e-7 &  3(3) &   9 & 1.3e-16 & 14(2) &  24 & 4.3e-15 &  9(3) &  25 & 1.2e-16 &  5(5) &  24   \\
         & 1.0e1 & 2.7e-6 &  6(2) &  12 & 3.0e-8 &  4(4) &  12 & 2.5e-15 & 20(2) &  33 & 2.9e-16 & 12(4) &  33 & 1.8e-15 &  5(5) &  33   \\
         & 3.2e1 & 2.7e-6 & 10(2) &  13 & 3.0e-7 &  4(4) &  13 & 8.2e-15 & 36(2) &  42 & 4.8e-15 & 12(4) &  44 & 7.3e-17 &  6(6) &  42   \\
         & 1.0e2 & 4.8e-6 & 18(2) &  16 & 2.7e-8 &  5(5) &  15 & 6.2e-15 & 74(2) &  49 & 2.1e-16 & 15(5) &  49 & 5.7e-17 &  6(6) &  49   \\
         & 3.2e2 & 2.4e-6 & 18(3) &  16 & 2.6e-8 &  5(5) &  16 & 9.9e-15 & 66(3) &  64 & 2.3e-16 & 15(5) &  56 & 9.4e-17 &  7(7) &  56   \\
         & 1.0e3 & 2.8e-4 & 21(3) & -- & 3.8e-8 &  5(5) &  16 & 1.6e-4 & 30(3) & -- & 2.6e-16 & 15(5) &  63 & 9.7e-17 &  7(7) &  62   \\
         & 3.2e3 & 1.4e-3 & 15(5) & -- & 3.6e-8 &  5(5) &  16 & 2.1e-3 & 15(5) & -- & 1.3e-15 & 15(5) &  69 & 1.1e-16 &  7(7) &  68   \\
\midrule
$1000$ & 1.3e0 & 2.5e-7 &  4(2) &   6 & 2.4e-7 &  2(2) &   6 & 5.2e-14 & 10(2) &  22 & 2.3e-16 &  6(2) &  12 & 1.9e-17 &  4(4) &  12   \\
         & 3.2e0 & 1.3e-6 &  4(2) &   6 & 1.3e-7 &  3(3) &   6 & 2.6e-14 & 12(2) &  32 & 1.8e-15 &  9(3) &  22 & 1.8e-17 &  5(5) &  22   \\
         & 1.0e1 & 1.2e-5 &  4(2) &   9 & 7.3e-6 &  3(3) &   9 & 7.9e-15 & 18(2) &  31 & 9.6e-15 & 12(3) &  31 & 8.6e-16 &  5(5) &  31   \\
         & 3.2e1 & 2.5e-5 &  6(2) &  11 & 1.2e-7 &  4(4) &  10 & 8.1e-14 & 32(2) &  93 & 1.9e-15 & 12(4) &  39 & 2.3e-17 &  6(6) &  39   \\
         & 1.0e2 & 3.7e-5 &  8(2) &  13 & 5.7e-7 &  4(4) &  11 & 7.9e-14 & 72(2) &  48 & 1.3e-16 & 16(4) &  48 & 2.9e-17 &  6(6) &  48   \\
         & 3.2e2 & 2.4e-4 & 102(2) & -- & 7.4e-9 &  5(5) &  12 & 2.4e-4 & 102(2) & -- & 1.0e-16 & 15(5) &  55 & 2.7e-17 &  7(7) &  55   \\
         & 1.0e3 & 3.4e-4 & 21(3) & -- & 8.3e-9 &  5(5) &  12 & 5.5e-4 & 21(3) & -- & 1.9e-16 & 15(5) &  62 & 2.9e-17 &  7(7) &  61   \\
         & 3.2e3 & 5.5e-5 &  9(3) &  15 & 1.1e-8 &  5(5) &  10 & 2.2e-7 & 75(3) & -- & 1.3e-15 & 15(5) &  75 & 2.9e-17 &  7(7) &  69   \\
\bottomrule
\end{tabularx} 	}
	\label{tab:test-synthetic-cond-small}
\end{table}

The results are presented in Table~\ref{tab:test-synthetic-cond-small}.
Among successful runs, the iteration counts and retained ranks generally increase with the condition number, and the two formulations have broadly similar behaviour, though the $\ldlt$-type IR often requires fewer Newton steps when both formulations converge with solver precision $\us=$ bf16.
With $u_s=$ fp32, both algorithms converge in both working precisions for all the problems shown; additional tests (not presented) with $u=$ fp64 showed convergence for condition numbers up to about $10^7$.
In contrast, decreasing the solver precision $u_s$ to bf16 limits the range of problems over which the IR scheme converges.
For $n=100$, both algorithms converge for all tested condition numbers up to about $10^{2.5}$ in both working precisions.
For $n=1000$, the corresponding ranges extend to about $10^2$ for the $\ldlt$-type IR in both working precisions, and to about $10^{1.5}$ and $10^1$ for the Cholesky-type IR in single and double working precisions, respectively.
Both algorithms also converge for a problem with condition number approximately $10^{3.5}$ in single working precision, despite failing at the intermediate tested condition numbers $10^{2.5}$ and $10^3$.

The results are largely in agreement with the condition number bounds in Table~\ref{table:ir3}, but they also display the instability of the sign function Newton iteration in floating-point arithmetic, which occurs when the matrix $A$ has eigenvalues close to the imaginary axis~\cite{bhm97} and may be indicated by a large $\kappa_{\sign}(B)$~\cite[sect.~5.1]{high08}, where
$B = \begin{bsmallmatrix}
	A & W \\ 0 &  -A^T
\end{bsmallmatrix}$.
In particular, we found that the solution update of Fragment~\ref{alg.sol-update-chol} and Fragment~\ref{alg.sol-update-ldlt} has been performed accurately with $u_c=u$, where the $b_3$ of~\eqref{eq:assump-c3-ldlt} remains moderate and approaches $1$ as the refinement proceeds.

\subsection{Conditional speedups for different cost ratios}\label{sect:expected-speedup}

\begin{figure}[!t]
	\centering
	\input{speedup_thresholds.tex}
	\caption{Relative residuals at termination (left) and conditional break-even cost ratios (right) for the $\ldlt$-IR experiments in Table~\ref{tab:test-synthetic-cond-small}. Threshold markers above a horizontal $\rho$ line predict a speedup at that cost ratio under the models of Section~\ref{sect:mp-ir-cost}.}
	\label{fig:synthetic-speedup-thresholds}
\end{figure}

Figure~\ref{fig:synthetic-speedup-thresholds} plots the attained relative residuals~\eqref{eq:res-fac-sol} and corresponding break-even cost ratios from the $\ldlt$-IR runs in Table~\ref{tab:test-synthetic-cond-small}, indicating when using bf16 for the inner solver is predicted to reduce the time spent on the dominant computational work.
The upper row compares solver precisions $u_s=$ bf16 and fp32 at fixed working precision $u=$ fp32; the lower row compares $u_s=$ bf16 and fp64 at $u=$ fp64.
Accordingly, $\rho$ in~\eqref{eq:rho-cost} denotes the bf16-to-fp32 cost ratio in the upper row and the bf16-to-fp64 cost ratio in the lower row.
Each panel shows the $n=100$ results in its left half and the $n=1000$ results in its right half, with a separate logarithmic scale for the condition number in each half.
For each successful precision pair, the markers on the right give $\rho_{\sigma}^{*}$ and $\rho_{\max}^{*}$ from~\eqref{eq:speedup-no-reuse} and~\eqref{eq:speedup-reuse}, using total and maximal inner iteration counts, respectively.
The three horizontal lines mark the hypothetical cost ratios $\rho=1/2$, $\rho=1/4$, and $\rho=1/16$, respectively.
Crosses identify bf16 runs that terminate without attaining the prescribed residual tolerance; those at the bottom of the right panels have no numerical threshold value and do not have a speedup estimate.

Markers above the horizontal $\rho$ lines in the right panels of Figure~\ref{fig:synthetic-speedup-thresholds} illustrate how a smaller $\rho$ can offset the additional Newton steps required by a lower-precision solver.
For $n=1000$ and $\cond=1.3$, bf16 versus fp32 solver precision at $u=$ fp32 moves from break-even to speedups of $2$ and $8$ without reuse ($2$, $4$, and $16$ with reuse) as $\rho$ decreases from $1/2$ to $1/4$ and $1/16$.
At $u=$ fp64, bf16 versus fp64 instead moves from an estimated speedup ratio of $0.8$, indicating a slowdown, to speedups of $1.6$ and $6.4$ ($4$, $8$, and $32$ with reuse).
Similarly, the fp32/fp64 solver comparison in Table~\ref{tab:test-synthetic-cond-small} for $n=1000$ and $\cond=10^2$ at $u=$ fp64 becomes beneficial as $\rho$ decreases, with the predicted speedup ratio increasing from $0.75$ to $1.5$ and $6$ ($3$, $6$, and $24$ with reuse), for the same three values of $\rho$.
Nevertheless, the estimated speedups for the reuse case require caution if $\vartheta k^{\sigma}\ll k^{\max}$ fails for either run (cf.~Section~\ref{sect:mp-ir-cost-alt}).
For example, the bf16/fp64 $\ldlt$-IR comparison with $n=1000$, $\cond=10^2$, and $u=$ fp64 in Table~\ref{tab:test-synthetic-cond-small} yields reuse estimates of $6$, $12$, and $48$ at $\rho=1/2$, $1/4$, and $1/16$, respectively.
Despite these large predicted speedups, the bf16 run has $\vartheta k^{\sigma}/k^{\max}=0.1\times72/2=3.6$, so the condition is not satisfied.
In this case, the cumulative cost of rank truncations and matrix multiplications used to update the low-rank factors may then be significant, so the estimates may cease to be a good indicator for the overall speedup.

\subsection{Accuracy on benchmark problems}

Finally, we evaluate the accuracy of the IR algorithms on Lyapunov equations from the SLICOT library\footnote{\url{https://www.slicot.org}} of benchmark examples of model reduction problems~\cite{chva02}.
The experiments here and in the following subsection reset the pseudorandom number generator with \texttt{rng(0,'twister')} before each complete run of IR.
Key characteristics of the test problems are listed in Table~\ref{tab:testmat-slicot}.
For each dataset, we estimate $\kappa_F(A)$ as well as $\kappa_{\sign}(B)$ in the Frobenius norm.
The sign-function condition numbers are estimated in fp64 using six power iterations in \texttt{funm\_condest\_fro} from the Matrix Function Toolbox~\cite[App.~D]{high08}, supplied with analytic Fr\'echet derivative and adjoint actions evaluated through Sylvester equations.
Since the sign function Newton iteration~\eqref{eq:newton-lyapu} for solving the Lyapunov equation~\eqref{eq:lyap} is essentially computing $\sign(B)$, $\kappa_{\sign}(B)$ provides an indicator of its unstructured sensitivity. The convergence of the IR framework also depends on the residual accuracy attained by the inner solver.

\begin{table}[!htb]\caption{Summary of the test Lyapunov equations from the SLICOT library.}\label{tab:testmat-slicot}
	\renewcommand{\tabcolsep}{8pt}
	\centering
\begin{tabular}{lrrrll}
			\toprule
			Dataset  & $n$ & $m$ & Nonzeros & $\kappa_F(A)$\tnote{*}  & $\kappa_{\sign}(B)$\tnote{*} \\
			\midrule
			\texttt{beam}  	&   348	& 1 &  60,726 & $1.2\eu {7}$ & $1.4\eu {9}$ \\
			\texttt{build} 	&  	48  & 1	&  1,176 & $7.5\eu {4}$ & $1.9\eu {6}$ \\
			\texttt{CDplayer}&	120	 & 2 &  240   & $1.4\eu {5}$ & $4.1\eu {12}$ \\
			\texttt{eady} 	&	598 & 1 &  357,406 & $3.6\eu {3}$& $3.8\eu {6}$ \\
			\texttt{fom} 	& 	1,006 & 1 & 1,012 & $2.3\eu {4}$ & $5.2\eu {5}$ \\
			\texttt{heat-cont} & 200	& 1 &  598	  & $1.5\eu {5}$ & $1.0\eu {4}$ \\
			\texttt{iss} 	& 	270 & 3 &  405   & $2.5\eu {5}$ & $7.2\eu {6}$ \\
			\texttt{pde} 	& 	84  & 1	&  382   & $1.2\eu {2}$ & $1.6\eu {2}$ \\
			\texttt{random}  & 	200 & 1	&  2132  & $3.2\eu {3}$ & $1.2\eu {15}$ \\
			\bottomrule
		\end{tabular}
\end{table}

\begin{table}[!htb]\renewcommand{\tabcolsep}{4pt}
\centering
	\caption{Results on the problems presented in Table~\ref{tab:testmat-slicot}. A dash ``--'' in the rank indicates that the internal residual stopping criterion was not met.}
	\resizebox{1\linewidth}{!}{
		\input{table_slicot.tex} }\label{tab:test-slicot}
\end{table}

The numerical results are presented in Table~\ref{tab:test-slicot}, where the residuals reported for nonconvergent runs are computed from the factors returned at termination.
Perhaps not surprisingly, the IR framework with both types of solvers only reached convergence on few problems that are relatively well conditioned when $u_s=$ bf16. For the other problems, the limiting residual presented in Table~\ref{table:ir3} is clearly irrelevant.
With $u_s=$ bf16, each formulation converges on three of the nine problems at $u=$ fp32 and only on \texttt{pde} at $u=$ fp64.
With $u_s=$ fp32 and $u=$ fp64, the algorithm converges in most cases, with failures on \texttt{beam}, \texttt{iss}, and \texttt{random} for the Cholesky formulation, and on \texttt{eady} and \texttt{iss} for the $\ldlt$ formulation.
Since most problems within the dataset are mild- to ill-conditioned,
the total number of Newton iterations across all refinement steps is typically more than doubled when the solver precision $u_s$ decreases from fp64 to fp32. However, $k_{\mathrm{low}}^{\max}/k_{\mathrm{high}}^{\max}<1$ whenever both solvers converge at the same working precision, so the reuse model~\eqref{eq:speedup-reuse}, with $\rho\le1$, suggests potential acceleration on well-conditioned problems by using a lower solver precision.

\subsection{Using adaptive solver precision}\label{sect:adaptive-solver-precision}
The failures associated with a fixed low solver precision in the experiments motivate an adaptive precision-switching scheme for Algorithm~\ref{alg.ir3-lyap-ldlt}.
Within the scope of this work, we can start the inner solver in bf16 and switch its arithmetic to fp32, and subsequently to fp64 if needed, while keeping the working precision and precisions of residual factorization and solution update fixed at fp64.
The solver switches to the next higher precision if the inner Newton iteration terminates, either at the iteration limit or after stagnation of the inner coefficient iterates (see~\eqref{eq:newton-inner-stagnation}), without satisfying~\eqref{eq:newton-stop-tau}, evaluated using the current unit roundoff $\us$.
It also switches precision if two consecutive outer updates satisfy $\theta_i>0.9$ in~\eqref{eq:it-early-terminate} while the residual remains above the convergence tolerance $\tau_I$ specified in Section~\ref{sect:algorithmic-parameters}; with fixed solver precision, this condition terminates the refinement process.
At each precision switch, we retain the current solution factors and use their residual factors, computed in fp64 using Fragment~\ref{alg.res-fac-ldlt}, to solve~\eqref{eq:lyap-cor-ldlt} at the higher solver precision.
If either precision-switching condition is met at fp64, the algorithm reports nonconvergence as no higher solver precision is available.
The overall algorithm with adaptive precision switching is presented as Algorithm~\ref{alg:adaptive-precision-switching}.

\begin{algorithm2e}[t]
 \caption{$\ldlt$-type IR with adaptive precision switching.}\label{alg:adaptive-precision-switching}
 \SetAlgoVlined
 \Parameter{Relative convergence tolerance $\tau_I>0$, maximal number of corrections $i_{\max}\in\mathbb N^+$, and precisions $u=u_r=u_c=$ fp64}
 \KwIn{$A\in\R^{n\times n}$, $L\in\R^{n\times m}$, and symmetric PSD $S\in\R^{m\times m}$}
 \KwOut{Solution factors $Z$, $Y$ such that $X\approx ZYZ^T$}
 Set solver precision $u_s=$ bf16, $i=1$, and stagnation counter $c=0$\;
 Solve $AX+XA^T+LSL^T=0$ at precision $u_s$ and store $Z_1,Y_1$ at precision $u$\;
 $Z=Z_1$, $Y=Y_1$\;
 Evaluate $[L_1^\Delta,S_1^\Delta]=\textsc{ResFacLdlt}(A,L,S,Z_1,Y_1)$ using Fragment~\ref{alg.res-fac-ldlt} at precision $u_r$\;
 Set $\res_1=\operatorname{res}(Z_1,Y_1)$\;
 \While{$i\le i_{\max}$}{
  \If{$\res_i\le\tau_I$}{\textbf{break}\;}
  \If{the preceding inner solve fails~\eqref{eq:newton-stop-tau} or $c=2$}{
   \lIf{$u_s=$ fp64}{report nonconvergence and \textbf{break}}
   Switch $u_s$ to the next higher precision and set $c=0$\;
  }
  Round copies of $L_i^\Delta,S_i^\Delta$ to precision $u_s$ for the correction solve\;
  Solve~\eqref{eq:lyap-cor-ldlt} at precision $u_s$ and store $Z_i^\Delta,Y_i^\Delta$ at precision $u$\;
  \lIf{\eqref{eq:newton-stop-tau} is not satisfied}{discard $Z_i^\Delta,Y_i^\Delta$ and \textbf{continue}}
  Evaluate $[Z_{i+1},Y_{i+1}]=\textsc{SolUptLdlt}(Z_i,Y_i,Z_i^\Delta,Y_i^\Delta)$ using Fragment~\ref{alg.sol-update-ldlt} at precision $u_c$\;
  Evaluate $[L_{i+1}^\Delta,S_{i+1}^\Delta]=\textsc{ResFacLdlt}(A,L,S,Z_{i+1},Y_{i+1})$ using Fragment~\ref{alg.res-fac-ldlt} at precision $u_r$\;
  Set $\res_{i+1}=\operatorname{res}(Z_{i+1},Y_{i+1})$ and $\theta_{i+1}=\res_{i+1}/\res_i$\;
  \eIf{$\theta_{i+1}>0.9$}{$c=c+1$\;}{ $c=0$\;}
  $Z=Z_{i+1}$, $Y=Y_{i+1}$, $i=i+1$\;
 }
\end{algorithm2e}

We test Algorithm~\ref{alg:adaptive-precision-switching} on all SLICOT examples in Table~\ref{tab:testmat-slicot} except \texttt{pde} (which already converges with $u_s=$ bf16), using the $\ldlt$ formulation with $u=$ fp64 and initial solver precision $u_s=$ bf16.
We compare \emph{retaining} the current factors with \emph{restarting}: at each precision switch, the latter discards the factors and solves the original Lyapunov equation~\eqref{eq:lyap} afresh at the higher solver precision.
Newton steps are counted separately at each precision over all attempted inner solves.

Table~\ref{tab:slicot-adaptive-precision} compares the final residuals~\eqref{eq:res-fac-sol} and Newton iteration counts of the adaptive precision-switching scheme with those of the two fixed-solver-precision baselines.
Here $k_h^\sigma$, $k_s^\sigma$, and $k_d^\sigma$ denote the total numbers of Newton steps at solver precisions bf16, fp32, and fp64, respectively.
The fixed-fp32 solver attains $\tau_I$ on six of the eight problems, failing on \texttt{eady} and \texttt{iss}, while adaptive precision switching converges on seven when retaining the current factors and on all eight when restarting.
In particular, the seven successful runs with retained factors use no fp64 Newton steps.
For \texttt{eady}, the fixed-fp32 solver fails to converge, while the precision-switching algorithm with retained factors converges using only bf16 and fp32.
In this case, the iterates of the solution factors and the outer residual generated by the two algorithms followed different paths of convergence; the one with fixed-fp32 solver precision leads to stagnation, while the other happens to converge. For a similar reason, the adaptive precision-switching scheme converges on \texttt{iss} with restarting (after the final precision switching the algorithm has exactly the same behavior as $u = \us =$ fp64 in Table~\ref{tab:test-slicot}), while it fails to converge on the same problem with retained factors. This seems to suggest that increasing the solver precision cannot cure the divergence brought by a bad starting point at the Newton iterations.

For \texttt{iss}, the adaptive precision-switching scheme converges with restarting, while it fails to attain $\tau_I$ with retained factors even after switching to fp64.
This shows that increasing the solver precision alone does not always restore convergence when the current factors are retained.

\begin{table}[t]
 \centering
 \caption{Adaptive precision switching on the SLICOT examples with $u=$ fp64. A dagger marks a terminal residual exceeding the prescribed tolerance $\tau_I$.}
 \label{tab:slicot-adaptive-precision}
 \input{table_slicot_adaptive.tex}
\end{table}

Retaining the factors in the adaptive precision-switching scheme can in some cases reduce the subsequent higher-precision work.
For \texttt{beam} and \texttt{random}, it reduces the fp32 counts from $98$ to $84$ and from $75$ to $60$, respectively, compared with restarting or the fixed-fp32 solve.
Both adaptive variants perform $16$ bf16 steps on \texttt{beam} and $161$ on \texttt{random}.
On the other hand, this saving is not guaranteed, and it can also increase the subsequent refinement work.
On \texttt{build}, for example, retaining the factors requires $126$ fp32 steps compared with $70$ after a restart. 

\section{Conclusions}\label{sect.conclusions}
We have developed a mixed-precision IR framework for factored solutions of Lyapunov equations with low-rank right-hand sides, in Cholesky-type (Algorithm~\ref{alg.ir3-lyap-chol}) and $\ldlt$-type (Algorithm~\ref{alg.ir3-lyap-ldlt}) formulations.
The rounding error analysis is carried out for the general $\ldlt$ formulation, which includes the Cholesky formulation as a specialization.
It gives sufficient conditions for the attainable normwise residual and guides the choice of mixed precisions and other algorithmic parameters within the IR framework.
We then illustrate the IR framework using the sign function Newton iteration as the solver, and we developed a $\ldlt$-type formulation that enables the refinement of computed solution factors from an indefinite residual. 

The experiments on synthetic and model reduction benchmark problems assess the predicted accuracy and show the limitations of a fixed low solver precision as problem conditioning deteriorates.
Motivated by these limitations, we investigate an adaptive precision-switching scheme (Algorithm~\ref{alg:adaptive-precision-switching}) that increases the solver precision when stagnation is detected or an iteration limit is reached.
The results show that this scheme can replace some high-precision Newton steps with low-precision work, but further study of the scheme, including a more systematic strategy for precision switching, is beyond the scope of this work.

We also analyze the computational cost of the IR algorithms using the sign function Newton iteration as the solver and develop a flop-based model to estimate the potential speedup from using half precision in the solver.
The model accounts for the dominant computational overhead and the relative costs of low- and high-precision arithmetic.
Together with the rounding error analyses, this model can guide high-performance implementations of the mixed-precision IR framework on hardware that natively supports low-precision formats.
Developing these implementations and optimizing their performance using wall-clock measurements would require a separate study.

This work is the first step towards exploiting the emerging new reduced formats, such as the half precision, in solving the low-rank Lyapunov equation.
The \emph{solver-independent} error analysis also motivates using ADI-based~\cite{blp08} and Krylov-based methods~\cite{jbri06,simo07} as correction solvers within the mixed-precision IR framework.
A closely related direction is to develop mixed-precision variants of the compress-and-restart block Krylov methods of Kressner et~al.~\cite{klmp21}, whose Lyapunov variant uses residual factorization and solution update procedures similar to those of the present IR framework.
It may also be possible to accelerate the IR framework by using inexact Lyapunov solvers whose convergence tolerances are chosen adaptively as the refinement proceeds~\cite{kufr20}. 

\section*{Acknowledgments}\addcontentsline{toc}{section}{Acknowledgments}
The authors thank Massimiliano Fasi and Jonas Schulze for their helpful comments on a draft manuscript. 
The comments and suggestions from the anonymous referees helped improve the presentation of this paper.

\addcontentsline{toc}{section}{References}

\bibliographystyle{abbrvurl}
\bibliography{reference}

\end{document}

%% file: speedup_thresholds.tex
\begin{tikzpicture}
\begin{groupplot}[
 group style={group size=2 by 2,horizontal sep=1.85cm,vertical sep=3.3cm},
 width=0.43\linewidth,height=0.31\linewidth,
 xmode=log,log basis x=10,xmin=1,xmax=40000000,
 xtick={1,10,100,1000,10000,100000,1000000,10000000},
 xticklabels={$10^0$,$10^1$,$10^2$,$10^3$,$10^0$,$10^1$,$10^2$,$10^3$},
 xlabel={Condition number},xticklabel style={yshift=-2pt},clip=false,
 ymode=log,log basis y=10,grid=major,
 tick label style={font=\scriptsize},label style={font=\small},
 title style={font=\small,yshift=13pt},legend style={font=\scriptsize,draw=none,fill=none,at={(0.5,-0.32)},anchor=north,legend columns=2,column sep=5pt},
 legend cell align=left]
\nextgroupplot[title={bf16/fp32: residuals},ylabel={Relative residual},ymin=1e-09,ymax=0.003]
\node[font=\scriptsize,anchor=south] at (rel axis cs:0.237,1.025) {$n=100$};
\node[font=\scriptsize,anchor=south] at (rel axis cs:0.763,1.025) {$n=1000$};
\draw[gray!60,densely dashed] (rel axis cs:0.5,0) -- (rel axis cs:0.5,1);
\draw[white,line width=3pt] ([xshift=-3pt]rel axis cs:0.5,0) -- ([xshift=3pt]rel axis cs:0.5,0);
\draw[black,line width=0.4pt] ([xshift=-3pt,yshift=-2pt]rel axis cs:0.5,0) -- ([xshift=-1pt,yshift=2pt]rel axis cs:0.5,0);
\draw[black,line width=0.4pt] ([xshift=1pt,yshift=-2pt]rel axis cs:0.5,0) -- ([xshift=3pt,yshift=2pt]rel axis cs:0.5,0);
\draw[white,line width=3pt] ([xshift=-3pt]rel axis cs:0.5,1) -- ([xshift=3pt]rel axis cs:0.5,1);
\draw[black,line width=0.4pt] ([xshift=-3pt,yshift=-2pt]rel axis cs:0.5,1) -- ([xshift=-1pt,yshift=2pt]rel axis cs:0.5,1);
\draw[black,line width=0.4pt] ([xshift=1pt,yshift=-2pt]rel axis cs:0.5,1) -- ([xshift=3pt,yshift=2pt]rel axis cs:0.5,1);
\addplot+[black,dashed,no marks] coordinates {(1,5.96046447754e-06) (4000,5.96046447754e-06)};
\addlegendentry{$nu$}
\addplot+[black,dashed,no marks,forget plot] coordinates {(10000,5.96046447754e-05) (40000000,5.96046447754e-05)};
\addplot+[blue,only marks,mark=*,mark size=2pt,mark options={draw=blue,fill=white}] coordinates {(1.3,6.7e-07) (3.2,3.5e-07) (10,3e-08) (32,3e-07) (100,2.7e-08) (320,2.6e-08) (1000,3.8e-08) (3200,3.6e-08)};
\addlegendentry{fp32}
\addplot+[blue,only marks,mark=*,mark size=2pt,mark options={draw=blue,fill=white},forget plot] coordinates {(13000,2.4e-07) (32000,1.3e-07) (100000,7.3e-06) (320000,1.2e-07) (1000000,5.7e-07) (3200000,7.4e-09) (10000000,8.3e-09) (32000000,1.1e-08)};
\addplot+[black,only marks,mark=triangle*,mark size=2pt,mark options={draw=black,fill=white}] coordinates {(1.3,1.6e-06) (3.2,3.9e-06) (10,2.7e-06) (32,2.7e-06) (100,4.8e-06) (320,2.4e-06)};
\addlegendentry{bf16}
\addplot+[black,only marks,mark=triangle*,mark size=2pt,mark options={draw=black,fill=white},forget plot] coordinates {(13000,2.5e-07) (32000,1.3e-06) (100000,1.2e-05) (320000,2.5e-05) (1000000,3.7e-05) (32000000,5.5e-05)};
\addplot+[red!75!black,only marks,mark=x,mark size=2.5pt] coordinates {(1000,0.00028) (3200,0.0014)};
\addlegendentry{bf16 failure}
\addplot+[red!75!black,only marks,mark=x,mark size=2.5pt,forget plot] coordinates {(3200000,0.00024) (10000000,0.00034)};
\nextgroupplot[title={bf16/fp32: break-even ratios},ylabel={Cost ratio threshold $\rho^*$},ymin=0.01,ymax=4]
\node[font=\scriptsize,anchor=south] at (rel axis cs:0.237,1.025) {$n=100$};
\node[font=\scriptsize,anchor=south] at (rel axis cs:0.763,1.025) {$n=1000$};
\draw[gray!60,densely dashed] (rel axis cs:0.5,0) -- (rel axis cs:0.5,1);
\draw[white,line width=3pt] ([xshift=-3pt]rel axis cs:0.5,0) -- ([xshift=3pt]rel axis cs:0.5,0);
\draw[black,line width=0.4pt] ([xshift=-3pt,yshift=-2pt]rel axis cs:0.5,0) -- ([xshift=-1pt,yshift=2pt]rel axis cs:0.5,0);
\draw[black,line width=0.4pt] ([xshift=1pt,yshift=-2pt]rel axis cs:0.5,0) -- ([xshift=3pt,yshift=2pt]rel axis cs:0.5,0);
\draw[white,line width=3pt] ([xshift=-3pt]rel axis cs:0.5,1) -- ([xshift=3pt]rel axis cs:0.5,1);
\draw[black,line width=0.4pt] ([xshift=-3pt,yshift=-2pt]rel axis cs:0.5,1) -- ([xshift=-1pt,yshift=2pt]rel axis cs:0.5,1);
\draw[black,line width=0.4pt] ([xshift=1pt,yshift=-2pt]rel axis cs:0.5,1) -- ([xshift=3pt,yshift=2pt]rel axis cs:0.5,1);
\addplot+[blue,only marks,mark=*,mark size=2pt,mark options={draw=blue,fill=white}] coordinates {(1.3,0.5) (3.2,0.75) (10,0.666666666667) (32,0.4) (100,0.277777777778) (320,0.277777777778)};
\addlegendentry{No reuse}
\addplot+[blue,only marks,mark=*,mark size=2pt,mark options={draw=blue,fill=white},forget plot] coordinates {(13000,0.5) (32000,0.75) (100000,0.75) (320000,0.666666666667) (1000000,0.5) (32000000,0.555555555556)};
\addplot+[teal!65!black,only marks,mark=square*,mark size=2pt,mark options={draw=teal!65!black,fill=white}] coordinates {(1.3,1) (3.2,1.5) (10,2) (32,2) (100,2.5) (320,1.66666666667)};
\addlegendentry{Reuse}
\addplot+[teal!65!black,only marks,mark=square*,mark size=2pt,mark options={draw=teal!65!black,fill=white},forget plot] coordinates {(13000,1) (32000,1.5) (100000,1.5) (320000,2) (1000000,2) (32000000,1.66666666667)};
\addplot+[red!75!black,only marks,mark=x,mark size=2.5pt] coordinates {(1000,0.014) (3200,0.014)};
\addlegendentry{bf16 failure}
\addplot+[red!75!black,only marks,mark=x,mark size=2.5pt,forget plot] coordinates {(3200000,0.014) (10000000,0.014)};
\addplot+[gray,dotted,no marks,forget plot] coordinates {(1,0.5) (4000,0.5)};
\addplot+[gray,dotted,no marks,forget plot] coordinates {(10000,0.5) (40000000,0.5)};
\addplot+[gray,dashed,no marks,forget plot] coordinates {(1,0.25) (4000,0.25)};
\addplot+[gray,dashed,no marks,forget plot] coordinates {(10000,0.25) (40000000,0.25)};
\addplot+[gray,dashdotted,no marks,forget plot] coordinates {(1,0.0625) (4000,0.0625)};
\addplot+[gray,dashdotted,no marks,forget plot] coordinates {(10000,0.0625) (40000000,0.0625)};
\node[anchor=west,font=\scriptsize,xshift=3pt] at (axis cs:40000000,0.5) {$\rho=1/2$};
\node[anchor=west,font=\scriptsize,xshift=3pt] at (axis cs:40000000,0.25) {$\rho=1/4$};
\node[anchor=west,font=\scriptsize,xshift=3pt] at (axis cs:40000000,0.0625) {$\rho=1/16$};
\nextgroupplot[title={bf16/fp64: residuals},ylabel={Relative residual},ymin=1e-18,ymax=0.003,ytick={1e-18,1e-14,1e-10,1e-6,1e-3}]
\node[font=\scriptsize,anchor=south] at (rel axis cs:0.237,1.025) {$n=100$};
\node[font=\scriptsize,anchor=south] at (rel axis cs:0.763,1.025) {$n=1000$};
\draw[gray!60,densely dashed] (rel axis cs:0.5,0) -- (rel axis cs:0.5,1);
\draw[white,line width=3pt] ([xshift=-3pt]rel axis cs:0.5,0) -- ([xshift=3pt]rel axis cs:0.5,0);
\draw[black,line width=0.4pt] ([xshift=-3pt,yshift=-2pt]rel axis cs:0.5,0) -- ([xshift=-1pt,yshift=2pt]rel axis cs:0.5,0);
\draw[black,line width=0.4pt] ([xshift=1pt,yshift=-2pt]rel axis cs:0.5,0) -- ([xshift=3pt,yshift=2pt]rel axis cs:0.5,0);
\draw[white,line width=3pt] ([xshift=-3pt]rel axis cs:0.5,1) -- ([xshift=3pt]rel axis cs:0.5,1);
\draw[black,line width=0.4pt] ([xshift=-3pt,yshift=-2pt]rel axis cs:0.5,1) -- ([xshift=-1pt,yshift=2pt]rel axis cs:0.5,1);
\draw[black,line width=0.4pt] ([xshift=1pt,yshift=-2pt]rel axis cs:0.5,1) -- ([xshift=3pt,yshift=2pt]rel axis cs:0.5,1);
\addplot+[black,dashed,no marks] coordinates {(1,1.11022302463e-14) (4000,1.11022302463e-14)};
\addlegendentry{$nu$}
\addplot+[black,dashed,no marks,forget plot] coordinates {(10000,1.11022302463e-13) (40000000,1.11022302463e-13)};
\addplot+[blue,only marks,mark=*,mark size=2pt,mark options={draw=blue,fill=white}] coordinates {(1.3,1.5e-16) (3.2,1.2e-16) (10,1.8e-15) (32,7.3e-17) (100,5.7e-17) (320,9.4e-17) (1000,9.7e-17) (3200,1.1e-16)};
\addlegendentry{fp64}
\addplot+[blue,only marks,mark=*,mark size=2pt,mark options={draw=blue,fill=white},forget plot] coordinates {(13000,1.9e-17) (32000,1.8e-17) (100000,8.6e-16) (320000,2.3e-17) (1000000,2.9e-17) (3200000,2.7e-17) (10000000,2.9e-17) (32000000,2.9e-17)};
\addplot+[black,only marks,mark=triangle*,mark size=2pt,mark options={draw=black,fill=white}] coordinates {(1.3,1.2e-15) (3.2,1.3e-16) (10,2.5e-15) (32,8.2e-15) (100,6.2e-15) (320,9.9e-15)};
\addlegendentry{bf16}
\addplot+[black,only marks,mark=triangle*,mark size=2pt,mark options={draw=black,fill=white},forget plot] coordinates {(13000,5.2e-14) (32000,2.6e-14) (100000,7.9e-15) (320000,8.1e-14) (1000000,7.9e-14)};
\addplot+[red!75!black,only marks,mark=x,mark size=2.5pt] coordinates {(1000,0.00016) (3200,0.0021)};
\addlegendentry{bf16 failure}
\addplot+[red!75!black,only marks,mark=x,mark size=2.5pt,forget plot] coordinates {(3200000,0.00024) (10000000,0.00055) (32000000,2.2e-07)};
\nextgroupplot[title={bf16/fp64: break-even ratios},ylabel={Cost ratio threshold $\rho^*$},ymin=0.01,ymax=4]
\node[font=\scriptsize,anchor=south] at (rel axis cs:0.237,1.025) {$n=100$};
\node[font=\scriptsize,anchor=south] at (rel axis cs:0.763,1.025) {$n=1000$};
\draw[gray!60,densely dashed] (rel axis cs:0.5,0) -- (rel axis cs:0.5,1);
\draw[white,line width=3pt] ([xshift=-3pt]rel axis cs:0.5,0) -- ([xshift=3pt]rel axis cs:0.5,0);
\draw[black,line width=0.4pt] ([xshift=-3pt,yshift=-2pt]rel axis cs:0.5,0) -- ([xshift=-1pt,yshift=2pt]rel axis cs:0.5,0);
\draw[black,line width=0.4pt] ([xshift=1pt,yshift=-2pt]rel axis cs:0.5,0) -- ([xshift=3pt,yshift=2pt]rel axis cs:0.5,0);
\draw[white,line width=3pt] ([xshift=-3pt]rel axis cs:0.5,1) -- ([xshift=3pt]rel axis cs:0.5,1);
\draw[black,line width=0.4pt] ([xshift=-3pt,yshift=-2pt]rel axis cs:0.5,1) -- ([xshift=-1pt,yshift=2pt]rel axis cs:0.5,1);
\draw[black,line width=0.4pt] ([xshift=1pt,yshift=-2pt]rel axis cs:0.5,1) -- ([xshift=3pt,yshift=2pt]rel axis cs:0.5,1);
\addplot+[blue,only marks,mark=*,mark size=2pt,mark options={draw=blue,fill=white}] coordinates {(1.3,0.333333333333) (3.2,0.357142857143) (10,0.25) (32,0.166666666667) (100,0.0810810810811) (320,0.106060606061)};
\addlegendentry{No reuse}
\addplot+[blue,only marks,mark=*,mark size=2pt,mark options={draw=blue,fill=white},forget plot] coordinates {(13000,0.4) (32000,0.416666666667) (100000,0.277777777778) (320000,0.1875) (1000000,0.0833333333333)};
\addplot+[teal!65!black,only marks,mark=square*,mark size=2pt,mark options={draw=teal!65!black,fill=white}] coordinates {(1.3,2) (3.2,2.5) (10,2.5) (32,3) (100,3) (320,2.33333333333)};
\addlegendentry{Reuse}
\addplot+[teal!65!black,only marks,mark=square*,mark size=2pt,mark options={draw=teal!65!black,fill=white},forget plot] coordinates {(13000,2) (32000,2.5) (100000,2.5) (320000,3) (1000000,3)};
\addplot+[red!75!black,only marks,mark=x,mark size=2.5pt] coordinates {(1000,0.014) (3200,0.014)};
\addlegendentry{bf16 failure}
\addplot+[red!75!black,only marks,mark=x,mark size=2.5pt,forget plot] coordinates {(3200000,0.014) (10000000,0.014) (32000000,0.014)};
\addplot+[gray,dotted,no marks,forget plot] coordinates {(1,0.5) (4000,0.5)};
\addplot+[gray,dotted,no marks,forget plot] coordinates {(10000,0.5) (40000000,0.5)};
\addplot+[gray,dashed,no marks,forget plot] coordinates {(1,0.25) (4000,0.25)};
\addplot+[gray,dashed,no marks,forget plot] coordinates {(10000,0.25) (40000000,0.25)};
\addplot+[gray,dashdotted,no marks,forget plot] coordinates {(1,0.0625) (4000,0.0625)};
\addplot+[gray,dashdotted,no marks,forget plot] coordinates {(10000,0.0625) (40000000,0.0625)};
\node[anchor=west,font=\scriptsize,xshift=3pt] at (axis cs:40000000,0.5) {$\rho=1/2$};
\node[anchor=west,font=\scriptsize,xshift=3pt] at (axis cs:40000000,0.25) {$\rho=1/4$};
\node[anchor=west,font=\scriptsize,xshift=3pt] at (axis cs:40000000,0.0625) {$\rho=1/16$};
\end{groupplot}
\end{tikzpicture}

%% file: table_slicot.tex
\begin{tabularx}{1.37\textwidth}{@{\extracolsep{\fill}}l|lrr|lrr|lrr|lrr|lrr}
\toprule
 Cholesky-type & \multicolumn{6}{c|}{$u=$ fp32} & \multicolumn{9}{c}{$u=$ fp64} \\
\midrule
 & \multicolumn{3}{c|}{$u_s=$ bf16} & \multicolumn{3}{c|}{$u_s=$ fp32} & \multicolumn{3}{c|}{$u_s=$ bf16} & \multicolumn{3}{c|}{$u_s=$ fp32} & \multicolumn{3}{c}{$u_s=$ fp64} \\
\text{Dataset} & $\res$ & $\iter$ & $\rk$ & $\res$ & $\iter$ & $\rk$ & $\res$ & $\iter$ & $\rk$ & $\res$ & $\iter$ & $\rk$ & $\res$ & $\iter$ & $\rk$ \\
\midrule
\texttt{beam} & 8.1e-4 & 80(16) & -- & 2.8e-7 & 14(14) &  54 & 1.6e-2 & 48(16) & -- & 6.5e-6 & 42(14) & -- & 8.2e-16 & 16(16) & 134   \\
\texttt{build} & 1.5e-3 & 36(12) & -- & 7.9e-8 & 14(14) &  35 & 1.4e-3 & 36(12) & -- & 4.3e-15 & 84(14) &  48 & 2.5e-17 & 15(15) &  48   \\
\texttt{CDplayer} & 4.3e-6 &  8(8) &   2 & 1.4e-9 & 16(16) &  10 & 5.9e-8 & 112(8) & -- & 1.3e-16 & 48(16) & 116 & 1.9e-16 & 18(18) & 116   \\
\texttt{eady} & 2.9e-3 & 48(12) & -- & 1.6e-7 & 15(15) &  12 & 2.5e-3 & 48(12) & -- & 3.5e-14 & 45(15) & 130 & 2.4e-16 & 16(16) &  88   \\
\texttt{fom} & 3.9e-3 & 12(3) & -- & 6.5e-9 & 13(13) &  12 & 3.9e-3 & 12(3) & -- & 9.9e-16 & 39(13) &  31 & 7.7e-17 & 15(15) &  27   \\
\texttt{heat-cont} & 2.0e-2 & 12(4) & -- & 2.5e-8 &  7(7) &  10 & 2.3e-2 & 12(4) & -- & 2.6e-15 & 28(7) &  29 & 4.4e-17 &  9(9) &  26   \\
\texttt{iss} & 6.5e-6 & 30(15) &  12 & 2.0e-8 & 21(21) &  46 & 4.6e-2 & 60(15) & -- & 5.7e-8 & 63(21) & -- & 4.9e-17 & 23(23) & 223   \\
\texttt{pde} & 7.4e-7 &  8(2) &   5 & 2.6e-8 &  4(4) &   5 & 9.0e-16 & 22(2) &  11 & 1.1e-16 & 12(4) &  11 & 9.8e-17 &  6(6) &  11   \\
\texttt{random} & 5.3e-2 & 20(7) & -- & 6.2e-8 & 15(15) &   2 & 1.7e-2 & 20(7) & -- & 3.8e-9 & 150(15) & -- & 1.3e-16 & 16(16) &  24   \\
\midrule
 $\ldlt$-type & \multicolumn{6}{c|}{$u=$ fp32} & \multicolumn{9}{c}{$u=$ fp64} \\
\midrule
 & \multicolumn{3}{c|}{$u_s=$ bf16} & \multicolumn{3}{c|}{$u_s=$ fp32} & \multicolumn{3}{c|}{$u_s=$ bf16} & \multicolumn{3}{c|}{$u_s=$ fp32} & \multicolumn{3}{c}{$u_s=$ fp64} \\
\text{Dataset} & $\res$ & $\iter$ & $\rk$ & $\res$ & $\iter$ & $\rk$ & $\res$ & $\iter$ & $\rk$ & $\res$ & $\iter$ & $\rk$ & $\res$ & $\iter$ & $\rk$ \\
\midrule
\texttt{beam} & 5.0e-3 & 64(16) & -- & 2.9e-7 & 14(14) &  54 & 1.6e-2 & 80(16) & -- & 1.4e-15 & 98(14) & 134 & 8.2e-16 & 16(16) & 134   \\
\texttt{build} & 9.2e-4 & 36(12) & -- & 5.9e-8 & 14(14) &  35 & 2.6e-3 & 60(12) & -- & 1.4e-16 & 70(14) &  48 & 6.4e-17 & 15(15) &  48   \\
\texttt{CDplayer} & 4.1e-6 &  8(8) &   2 & 3.4e-8 & 16(16) &  10 & 3.4e-8 & 128(8) & -- & 7.9e-15 & 48(16) & 116 & 5.7e-17 & 18(18) & 116   \\
\texttt{eady} & 2.1e-5 & 72(12) &  21 & 1.3e-7 & 15(15) &  12 & 4.1e-7 & 204(12) & -- & 2.4e-9 & 75(15) & -- & 2.1e-16 & 16(16) &  88   \\
\texttt{fom} & 3.9e-3 & 12(3) & -- & 2.1e-7 & 13(13) &  12 & 3.9e-3 & 12(3) & -- & 2.3e-15 & 39(13) &  29 & 6.1e-16 & 15(15) &  27   \\
\texttt{heat-cont} & 1.5e-2 & 12(4) & -- & 3.5e-8 &  7(7) &  10 & 2.8e-2 & 12(4) & -- & 2.6e-16 & 28(7) &  26 & 6.2e-17 &  9(9) &  26   \\
\texttt{iss} & 3.1e-5 & 60(15) & -- & 1.8e-8 & 21(21) &  46 & 3.1e-5 & 60(15) & -- & 9.2e-8 & 63(21) & -- & 2.6e-17 & 23(23) & 223   \\
\texttt{pde} & 3.3e-7 &  8(2) &   5 & 2.7e-8 &  4(4) &   5 & 3.1e-15 & 18(2) &  11 & 1.1e-16 & 12(4) &  11 & 5.2e-17 &  6(6) &  11   \\
\texttt{random} & 9.9e-5 & 59(7) & -- & 3.8e-8 & 15(15) &   2 & 4.4e-5 & 161(7) & -- & 2.6e-16 & 75(15) &  24 & 1.0e-16 & 16(16) &  24   \\
\bottomrule
\end{tabularx}

%% file: table_slicot_adaptive.tex
\begin{tabular}{lrrrrrrrr}
\toprule
 & \multicolumn{2}{c}{Fixed fp32} & \multicolumn{2}{c}{Fixed fp64} & \multicolumn{2}{c}{Adaptive: retain} & \multicolumn{2}{c}{Adaptive: restart} \\
\cmidrule(lr){2-3}\cmidrule(lr){4-5}\cmidrule(lr){6-7}\cmidrule(lr){8-9}
Dataset & $k_s^\sigma$ & $\res$ & $k_d^\sigma$ & $\res$ & $(k_h^\sigma,k_s^\sigma,k_d^\sigma)$ & $\res$ & $(k_h^\sigma,k_s^\sigma,k_d^\sigma)$ & $\res$ \\
\midrule
\texttt{beam} & $98$ & 1.4e-15 & $16$ & 8.2e-16 & $(16,84,0)$ & 2.4e-15 & $(16,98,0)$ & 1.4e-15 \\
\texttt{build} & $70$ & 1.4e-16 & $15$ & 6.4e-17 & $(60,126,0)$ & 1.9e-15 & $(60,70,0)$ & 1.4e-16 \\
\texttt{CDplayer} & $48$ & 7.9e-15 & $18$ & 5.7e-17 & $(128,48,0)$ & 3.7e-16 & $(128,48,0)$ & 7.9e-15 \\
\texttt{eady} & $75$ & 2.4e-9$^{\dagger}$ & $16$ & 2.1e-16 & $(204,45,0)$ & 1.4e-15 & $(204,75,16)$ & 2.1e-16 \\
\texttt{fom} & $39$ & 2.3e-15 & $15$ & 6.1e-16 & $(12,39,0)$ & 1.1e-15 & $(12,39,0)$ & 2.3e-15 \\
\texttt{heat-cont} & $28$ & 2.6e-16 & $9$ & 6.2e-17 & $(12,35,0)$ & 1.0e-15 & $(12,28,0)$ & 2.6e-16 \\
\texttt{iss} & $63$ & 9.2e-8$^{\dagger}$ & $23$ & 2.6e-17 & $(60,63,46)$ & 9.3e-7$^{\dagger}$ & $(60,63,23)$ & 2.6e-17 \\
\texttt{random} & $75$ & 2.6e-16 & $16$ & 1.0e-16 & $(161,60,0)$ & 1.8e-16 & $(161,75,0)$ & 2.6e-16 \\
\bottomrule
\end{tabular}